\documentclass[a4paper,12pt]{article}

\usepackage{amsthm}
\usepackage{amsmath,amssymb,latexsym,amsfonts,mathrsfs}

\theoremstyle{plain}
\newtheorem{Thm}{Theorem}[section]
\newtheorem{Lem}[Thm]{Lemma}
\newtheorem{Prop}[Thm]{Proposition}
\newtheorem{Cor}[Thm]{Corollary}
\newtheorem{Conj}[Thm]{Conjecture}
\theoremstyle{definition}
\newtheorem{Def}[Thm]{Definition}
\newtheorem{Rem}[Thm]{Remark}

\newcommand{\Proof}[2][Proof]{\begin{proof}[{#1}] #2 \end{proof}}

\makeatletter
\renewcommand\section{\@startsection {section}{1}{\z@}%
                                   {-3.5ex \@plus -1ex \@minus -.2ex}%
                                   {2.3ex \@plus.2ex}%
                                   {\normalfont\large\bf}}
\makeatother

\makeatletter
\renewcommand\subsection{\@startsection {subsection}{1}{\z@}%
                                   {-3.5ex \@plus -1ex \@minus -.2ex}%
                                   {2.3ex \@plus.2ex}%
                                   {\normalfont\normalsize\bf}}
\makeatother

\numberwithin{equation}{section}

\makeatletter
\@addtoreset{footnote}{page}
\makeatother

\renewcommand{\d}{{\rm d}} 
\newcommand{\eps}{\ensuremath{\varepsilon}}
\newcommand{\e}{{\rm e}} 
\newcommand{\law}{\stackrel{{\rm law}}{=}}
\newcommand{\claw}{\stackrel{{\rm law}}{\longrightarrow}}
\newcommand{\tend}[2]{\mathrel{\mathop{\longrightarrow}\limits^{#1}_{#2}}}
\renewcommand{\hat}{\widehat}
\renewcommand{\tilde}{\widetilde}
\renewcommand{\bar}{\overline}
\newcommand{\un}[1]{\underline{#1}}

\newcommand{\absol}[1]{\left| #1 \right|} 
\newcommand{\rbra}[1]{\!\left( #1 \right)} 
\newcommand{\cbra}[1]{\!\left\{ #1 \right\}} 
\newcommand{\sbra}[1]{\!\left[ #1 \right]} 

\newcommand{\bE}{\ensuremath{\mathbb{E}}}
\newcommand{\bN}{\ensuremath{\mathbb{N}}}
\newcommand{\bP}{\ensuremath{\mathbb{P}}}
\newcommand{\bQ}{\ensuremath{\mathbb{Q}}}
\newcommand{\bR}{\ensuremath{\mathbb{R}}}
\newcommand{\cB}{\ensuremath{\mathcal{B}}}
\newcommand{\cD}{\ensuremath{\mathcal{D}}}
\newcommand{\cE}{\ensuremath{\mathcal{E}}}
\newcommand{\cF}{\ensuremath{\mathcal{F}}}
\newcommand{\vn}{\ensuremath{\mbox{{\boldmath $n$}}}}

\begin{document}

\begin{center}
{\large \bf 
Functional limit theorems for processes pieced together from excursions 
} 
\end{center}
\begin{center}
Kouji Yano\footnote{
Graduate School of Science, Kyoto University, JAPAN }\footnote{
The research of this author is supported by KAKENHI (26800058) 
and partially by KAKENHI (24540390).} 
\end{center}
\begin{center}
{\small \today}
\end{center}

\begin{abstract}
A notion of convergence of excursion measures is introduced. 
It is proved that convergence of excursion measures 
implies convergence in law of the processes pieced together from excursions. 
This result is applied to obtain homogenization theorems 
of jumping-in extensions 
for positive self-similar Markov processes, 
for Walsh diffusions 
and for the Brownian motion on the Sierpi\'nski gasket. 
\end{abstract}


\section{Introduction}

In the previous work \cite{MR2543582}, the author obtained homogenization results 
of jumping-in extensions for diffusion processes on the half line. 
The proof was based on the construction of a sample path from excursions 
using It\^o's excursion theory \cite{ItoKyoto} and the time-change method. 
The key to the proof was to prove convergence of time-changed paths of the Brownian excursion 
based on the results of Fitzsimmons--Yano \cite{MR2384481}. 

The aim of this paper is to establish a general limit theorem (Theorem \ref{glt}) 
which asserts, roughly speaking, that 
\begin{align}
\vn_n \to \vn_{\infty } \ \text{implies} \ X_n \claw X_{\infty } , 
\label{eq: imply}
\end{align}
where $ \vn_n $'s are excursion measures and 
$ X_n $'s are the processes pieced together from excursions. 
For a given Hunt process for which the origin is regular for itself, 
the excursion measure away from the origin characterizes the law of the Hunt process. 
Hence it may be natural that \eqref{eq: imply} should hold. 
But in what sense is ``$ \vn_n \to \vn_{\infty } $"? 

We introduce a notion of convergence of excursion measures 
as an analogue to Skorokhod's a.s.-convergence realization 
of weak convergence of probability measures. 
We roughly say that $ \vn_n \to \vn_{\infty } $ 
if all $ \vn_n $'s can be realized as the pullbacks of a common $ \sigma $-finite measure, 
say $ \vn_n = \nu \circ \Phi_n^{-1} $, 
where $ \Phi_n $'s are measurable mappings 
which take values in the functional space of c\`adl\`ag paths 
equipped with the Skorokhod topology 
and which satisfy $ \Phi_n \to \Phi_{\infty } $, $ \nu $-a.e. 
The key to the proof of \eqref{eq: imply} 
is to realize $ X_n $'s 
from a common Poisson point process 
and to construct random time-changes 
which establish the convergence $ \Phi_n \to \Phi_{\infty } $ in the Skorokhod topology. 

We apply the general theorem to obtain homogenization theorems of jumping-in extensions. 
Let $ S $ be a Borel subset of $ \bR^d $ containing 0 
and $ H^0 = \{ X,(\bP^0_x)_{x \in S} \} $ be a Hunt process stopped upon hitting 0. 
Let $ S' $ be a measurable space 
and let $ \{ \vn_v \}_{v \in S'} $ be a kernel 
such that for each $ v \in S' $ the measure $ \vn_v $ is the excursion measure 
of an extension of $ H^0 $. 
A jumping-in extension is the process $ X_{\rho,j} $ pieced together from excursions 
corresponding to the excursion measure defined by 
\begin{align}
\vn_{\rho,j} = \int_{S'} \rho(\d v) \vn_v + \int_{S \setminus \{ 0 \}} j(\d x) \bP^0_x 
\label{eq: vnrhoj}
\end{align}
for some finite measure $ \rho $ on $ S' $ 
and some $ \sigma $-finite measure $ j $ on $ S \setminus \{ 0 \} $. 
(The excursion measure of any extension of $ H^0 $ 
may admit a representation of the form \eqref{eq: vnrhoj}; 
see It\^o \cite[Section 7]{MR0402949}.) 
Let $ c>1 $ be a fixed constant. 
For $ \gamma>0 $, we define the scaling transformation 
\begin{align}
(\Psi_{\gamma } w)(t) =& c^{-\gamma } w(ct) . 
\label{eq: Psidef}
\end{align}
For certain constants $ \alpha >0 $ and $ \gamma >0 $, 
we study the following scaled objects: 
\begin{align}
\vn^{(n)}_{\rho,j} = c^{\gamma n} \vn_{\rho,j} \circ (\Psi_{\alpha }^n)^{-1} , 
\quad 
X^{(n)}_{\rho,j} = \Psi_{\alpha }^n X_{\rho,j} . 
\label{}
\end{align}
We shall provide sufficient conditions for the following two types of convergences: 
\begin{align}
\vn^{(n)}_{\rho,j} \to 
\begin{cases}
\vn_{\rho^*,0} & \text{in the {\bf jumping-in vanishing case}}, \\
\vn_{0,j^*} & \text{in the {\bf jumping-in dominant case}} 
\end{cases}
\label{eq: excm conv}
\end{align}
for some $ \rho^* $ and $ j^* $. 
Thanks to the general theorem \eqref{eq: imply}, 
the convergence \eqref{eq: excm conv} 
leads to 
the corresponding convergence in law of the scaled process $ X^{(n)}_{\rho,j} $, 
which can be regarded as a homogenization result. 
In particular, we take up positive self-similar Markov processes, 
Walsh diffusions, 
and the Brownian motion on the Sierpi\'nski gasket.

Let us give a remark on earlier works about jumping-in extensions. 
Jumping-in extensions of diffusion processes 
were discussed by Feller \cite{MR0092046} 
in his study of determination of all possible boundary conditions 
for the generator of a diffusion process with accessible boundaries. 
Such processes appear in the study of population genetics; 
see, e.g., Hutzenthaler--Taylor \cite{MR2796680}. 
The sample path construction of the jumping-in extensions 
was first established by It\^o--McKean \cite{MR0154338} 
for Brownian motions 
using time-change method involving an independent Poisson process. 
It\^o \cite{ItoKyoto} established his theory of Poisson point process of excursions 
to construct a sample path by piecing together from excursions 
produced by a Poisson point process. 
Yano utilized It\^o's method 
in \cite{MR2543582} to obtain homogenization results 
of jumping-in extensions for diffusion processes on the half line 
and in \cite{MR3192547} to determine possible jumping-in extensions 
of diffusion processes on an interval. 

Note also that Lambert--Simatos \cite{LS} proved \eqref{eq: imply} 
in a certain sense which is different from ours. 
They gave a general condition 
for convergence of regenerative processes 
assuming the convergence of excursions bigger than $ \eps $ in some given functional, 
which are called the {\em $ \eps $-big excursions}, for all $ \eps>0 $. 

This paper is organized as follows. 
In Section \ref{sec: glt}, we give basic facts 
about the piecing procedure of excursions 
and state the general limit theorem. 
In Section \ref{sec: hom}, we state the homogenization theorems 
for jumping-in extensions in a rather general framework. 
In Section \ref{sec: ex}, we discuss three examples of the homogenization theorems, 
the positive self-similar Markov processes, 
the Walsh diffusions and the Brownian motion on the Sierpi\'nski gasket. 
Sections \ref{sec: pr glt} and Section \ref{sec: pr hom} 
are devoted to the proofs of the general limit theorem 
and the homogenization theorems, respectively.

\section*{Acknowledgement} 
The author thanks Naotaka Kajino 
for his valuable suggestions on Lemma \ref{kajino}. 
The author also thanks the referee 
for pointing out an error of earlier versions.

\section{General limit theorem} \label{sec: glt}

\subsection{Notations about excursions}

Let $ d \ge 1 $ 
and let $ D=D_{\bR^d} $ denote the set of all functions $ w:[0,\infty ) \to \bR^d $ 
which are c\`adl\`ag, i.e., right continuous with left limits. 
We say $ w_n \to w $ {\em uniformly on compacts}, or simply $ w_n \to w $ {\em uc}, if 
$ \sup_{t \in [0,t_0]} |w_n(t)-w(t)| \to 0 $ for all $ t_0>0 $. 
We equip $ D $ with the {\em Skorokhod topology}; 
we say $ w_n \to w $ in $ D $ if 
there exists a sequence of time-changes $ \{ I_n \} $ of $ [0,\infty ) $ 
such that 
each $ I_n:[0,\infty ) \to [0,\infty ) $ is bijective, continuous and increasing 
and 
\begin{align}
I_n - I \to 0 \ \text{uc and} \ 
w_n - w \circ I_n \to 0 \ \text{uc}, 
\label{}
\end{align}
where $ I(t) \equiv t $ denotes the identity time-change. 
It is well-known that $ D $ is a Polish space. 
We write $ \cB(D) $ for the $ \sigma $-field generated by all open subsets of $ D $. 
Let $ X=(X(t))_{t \ge 0} $ denote the coordinate process on $ D $, i.e., 
\begin{align}
X(w)(t) = X(t)(w) = w(t) . 
\label{}
\end{align}
For $ x \in \bR^d $, we denote the hitting time of $ x $ by 
\begin{align}
T_x(w) = \inf \{ t > 0 : w(t) = x \} , 
\label{}
\end{align}
where we adopt the usual convention $ \inf \emptyset = \infty $. 
We denote 
\begin{align}
\| w \| = \sup_{t \ge 0} |w(t)| 
\quad \text{for $ w \in D $}. 
\label{}
\end{align}
Paths stopped upon hitting 0 are called {\em excursions away from $ 0 $}. 
The set of all excursions away from 0 will be denoted by 
\begin{align}
D^0 = \{ w \in D : w(t)=w(t \wedge T_0(w)) \ \text{for all $ t \ge 0 $} \} . 
\label{}
\end{align}
We write $ o \in D^0 $ for the path $ o(t) \equiv 0 $. 
Note that, for $ w \in D^0 $, we have $ T_0(w)=0 $ if and only if $ w=o $. 

For $ t \in [0,\infty ) $, we define the shift operator $ \theta_t:D \to D $ by 
\begin{align}
(\theta_t w)(s) = w(t+s) 
, \quad s \ge 0 . 
\label{}
\end{align}
We denote $ \cF^0_t = \sigma(X(s):s \le t) $ and set 
$ \cF_t = \bigcap_{\eps>0} \cF^0_{t+\eps} $.

\subsection{The process pieced together from excursions}

We denote $ \sharp \{ \cdot \} $ 
by the number of elements of the set $ \{ \cdot \} $. 
For a $ \sigma $-finite measure $ \nu $ on a measurable space $ E $ 
and a measurable functional $ f $ on $ E $, 
we write $ \nu[f] $ for $ \int_E f \d \nu $ whenever the integral is well-defined. 

We first recall the usual notion of a Poisson point process; 
see, e.g., \cite[\S I.9]{IW} for the basic facts about it. 
Let $ \nu $ be a $ \sigma $-finite measure on a measurable space $ E $. 
We call $ \{ (p^{(l)})_{l \in \cD(p)},\bP \} $ 
a {\em Poisson point process on $ E $ with characteristic measure $ \nu $} 
if the random measure $ N_l $ for $ l \ge 0 $ defined by 
\begin{align}
N_l(A) = \sharp \{ s \in \cD(p) \cap [0,l] : p^{(s)} \in A \} 
, \quad A \in \cB(E) 
\label{}
\end{align}
satisfies that 
for any non-negative measurable functional $ f $ on $ E $ 
the process $ (N_l[f])_{l \ge 0} $ is a Poisson process with intensity $ \nu[1-\e^{-f}] $. 

We second introduce an auxiliary notation modifying 
the usual notation of a Poisson point process. 
Let $ \vn $ be a $ \sigma $-finite measure on $ D $ such that $ \vn(\{ o \})=0 $. 
We call $ \{ p=(p^{(l)})_{l \ge 0},\bP \} $ 
a {\em Poisson point process on $ D $ outside $ o $ 
with characteristic measure $ \vn $} 
if $ \{ (p^{(l)})_{l \in \cD(p)},\bP \} $ for $ \cD(p) = \{ l \ge 0 : p^{(l)} \neq o \} $ 
is a Poisson point process on $ D \setminus \{ o \} $ 
with characteristic measure $ \vn|_{D \setminus \{ o \}} $. 
Note that a Poisson point process $ (p^{(l)})_{l \in \cD(p)} $ on $ D \setminus \{ o \} $ 
can always be extended to 
a Poisson point process on $ D $ outside $ o $ 
by putting $ p^{(l)}=o $ when $ l \notin \cD(p) $. 

Let $ (\vn,\varsigma) $ be the pair 
consisting of a $ \sigma $-finite measure $ \vn $ on $ D $ 
such that $ \vn(\{ o \}) = 0 $
and a non-negative constant $ \varsigma $. 
Let $ p=(p^{(l)})_{l \ge 0} $ be a Poisson point process on $ D $ outside $ o $ 
with characteristic measure $ \vn $. 
Noting that $ p^{(l)} \in D $ for all $ l \ge 0 $, 
we have 
\begin{align}
T_0(p^{(l)}) = \inf \{ t>0 : p^{(l)}(t)=0 \} . 
\label{}
\end{align}
Let $ \varsigma \ge 0 $ be a constant 
and for $ l \ge 0 $ we define 
\begin{align}
\eta(l) = \eta(p,\varsigma;l) = \varsigma l + \sum_{s \le l} T_0(p^{(s)}) . 
\label{eq: eta}
\end{align}
We introduce the following conditions on the pair $ (\vn,\varsigma) $: 
\begin{description}
\item[(N0)] 
$ X \in D^0 $ and $ 0<T_0<\infty $, $ \vn $-a.e.; 
\item[(N1)] 
$ \vn[T_0 \wedge 1] < \infty $; 
\item[(N2)] 
either $ \varsigma > 0 $ or $ \vn(D) = \infty $; 
\item[(N3)] 
$ \vn (\| X \| \ge r) < \infty $ for all $ r>0 $. 
\end{description}
If the conditions {\bf (N0)} and {\bf (N1)} are satisfied, 
we see that $ p^{(l)} \in D^0 $ for all $ l \ge 0 $ and that 
$ T_0(p^{(l)})=0 $ for all but countably many $ l $.

\begin{Lem} \label{lem: ilt}
Suppose that the conditions {\bf (N0)} and {\bf (N1)} are satisfied. 
Then $ (\eta(l))_{l \ge 0} = (\eta(p,\varsigma;l))_{l \ge 0} $ 
is an increasing L\'evy process with Laplace transform 
\begin{align}
\bE [\e^{- \lambda \eta(l)}] = 
\exp \cbra{ - l \varsigma \lambda - l \vn \sbra{1 - \e^{- \lambda T_0}} } 
, \quad \lambda \ge 0 . 
\label{}
\end{align}
If, moreover, the condition {\bf (N2)} is satisfied, then 
it is strictly increasing. 
\end{Lem}

Lemma \ref{lem: ilt} is well-known, and so we omit its proof. 

The following proposition enables us to piece a process together from excursions. 

\begin{Prop} \label{piece}
Suppose that the conditions {\bf (N0)}-{\bf (N2)} are satisfied. 
Define 
\begin{align}
L(t) = L(p,\varsigma;t) = \inf \{ l \ge 0 : \eta(l)>t \} 
, \quad t \ge 0 
\label{}
\end{align}
and 
\begin{align}
X(t) = X(p,\varsigma;t) = 
\begin{cases}
p^{(l)}(t-\eta(l-)) 
& \text{if $ \eta(l-) \le t < \eta(l) $ for some $ l \ge 0 $}, \\
0 & \text{otherwise}. 
\end{cases}
\label{}
\end{align}
Then it holds that 
\begin{align}
\int_0^t 1_{\{ X(s)=0 \}} \d s = \varsigma L(t) . 
\label{eq: occup}
\end{align}
If, moreover, the condition {\bf (N3)} is satisfied, 
then the process $ X(p,\varsigma)=(X(p,\varsigma;t))_{t \ge 0} $ is $ D $-valued. 
If $ \vn(D)=\infty $, the converse is also true: 
if the process $ X(p,\varsigma) $ is $ D $-valued 
then the condition {\bf (N3)} is satisfied. 
\end{Prop}

The proof of Proposition \ref{piece} will be given in Section \ref{sec: piece}.

\subsection{General limit theorem}

For real-valued measurable functions $ f_1,f_2,\ldots $ and $ f_{\infty } $ 
defined on a measure space $ (E,\cE,\nu) $, 
we say that $ f_n \to f_{\infty } $, $ \nu $-{\em almost uniformly} 
if for any $ \eps>0 $ there exists $ A \in \cE $ such that 
$ \nu(A) < \eps $ and $ \sup_{A^c} |f_n-f_{\infty }| \to 0 $. 
Imitating the Skorokhod representation of almost sure convergence, 
we introduce the following notion of convergence. 

\begin{Def} \label{def: conv}
Let $ \vn_1,\vn_2,\ldots $ and $ \vn_{\infty } $ 
be $ \sigma $-finite measures on $ D $. 
We say that $ \vn_n \to \vn_{\infty } $ 
if there exist a Polish space $ E $, 
a $ \sigma $-finite measure $ \nu $ on $ E $ 
and measurable mappings $ \Phi_1,\Phi_2,\ldots,\Phi_{\infty } $ 
from $ E $ to $ D $ such that the following conditions hold: 
\begin{description}
\item[(G1)] 
$ \vn_n = (\nu \circ \Phi_n^{-1})|_{D \setminus \{ o \}} $ 
for $ n=1,2,\ldots $ and $ \infty $; 
\item[(G2)] 
$ \Phi_n \to \Phi_{\infty } $ in $ D $, $ \nu $-a.e.; 
\item[(G3)] 
$ \| \Phi_n \| \to \| \Phi_{\infty } \| $, $ \nu $-almost uniformly; 
\item[(G4)] 
$ T_0(\Phi_n) \to T_0(\Phi_{\infty }) $, 
$ \nu $-a.e.; 
\item[(G5)] 
there exists $ N \in \bN $ such that 
$ \nu \sbra{ 1 \wedge \sup_{n \ge N} T_0(\Phi_n) } < \infty $. 
\end{description}
\end{Def}

We shall see in Lemma \ref{lem2} that 
Condition {\bf (G3)} can be replaced by the following: 
\begin{description}
\item[(G3)$ ' $] 
$ \nu \rbra{ \bigcup_{n=1}^{\infty } \cbra{ \| \Phi_n \| \ge r } } < \infty $ 
for all $ r>0 $. 
\end{description}

\begin{Rem}
Condition {\bf (G2)} does not imply Condition {\bf (G4)}. 
This is because the functional $ T_0:D \to [0,\infty ] $ is not continuous; 
for instance, 
the sequence of functions $ w_n \in D_{\bR} $ defined by 
\begin{align}
w_n(t) = 
\begin{cases}
1 & \text{if $ 0 \le t < 1 $}, \\
1/n & \text{if $ 1 \le t < 2 $}, \\
0 & \text{if $ t \ge 2 $}, 
\end{cases}
\label{}
\end{align}
satisfies that $ w_n $ converges to $ w_{\infty } $ in $ D_{\bR} $, 
while $ T_0(w_n) \equiv 2 \not\to 1 = T_0(w_{\infty }) $. 
\end{Rem}

\begin{Thm} \label{glt}
Let $ \vn_1,\vn_2,\ldots $ and $ \vn_{\infty } $
be $ \sigma $-finite measures on $ D $. 
Let $ \varsigma_1,\varsigma_2,\ldots $ and $ \varsigma_{\infty } $ be non-negative constants. 
Suppose that the following conditions hold: 
\begin{description}
\item[(A1)] 
for each $ n \in \bN \cup \{ \infty \} $, the pair $ (\vn_n,\varsigma_n) $ satisfies 
Conditions {\bf (N0)}-{\bf (N3)}; 
\item[(A2)] 
$ \vn_n \to \vn_{\infty } $; 
\item[(A3)] 
$ X(T_0-) = 0 $, $ \vn_{\infty } $-a.e.; 
\item[(A4)] 
$ \varsigma_n \to \varsigma_{\infty } $. 
\end{description}
For $ n \in \bN \cup \{ \infty \} $, 
let $ p_n $ be a Poisson point process on $ D $ outside $ o $ 
with characteristic measure $ \vn_n $. 
Denote $ \eta_n(l) = \eta(p_n,\varsigma_n;l) $ and $ X_n(t) = X(p_n,\varsigma_n;t) $. 
Then it holds that 
\begin{align}
(X_n,L_n,\eta_n) \claw (X_{\infty },L_{\infty },\eta_{\infty }) 
\quad \text{as $ n \to \infty $} , 
\label{}
\end{align}
where the convergence is in the sense of law on $ D \times D_{\bR} \times D_{\bR} $. 
\end{Thm}

The proof of Theorem \ref{glt} will be given in Section \ref{sec: prf glt}.

\section{Homogenization theorems} \label{sec: hom}

\subsection{Excursion measures}

Let $ S $ be a Borel subset of $ \bR^d $ containing 0 
and let $ H^0= \{ X,(\bP^0_x)_{x \in S} \} $ be a Hunt process stopped upon hitting 0. 
A Hunt process $ H = \{ X, (\bP_x)_{x \in S} \} $ is called 
an {\em extension} of $ H^0 $ 
if the law of the stopped process $ X(t \wedge T_0) $ under $ \bP_x $ 
coincides with $ \bP^0_x $ for all $ x \in S $. 
We introduce the following set of conditions 
for an extension $ H $ of $ H^0 $: 
\begin{description}
\item[(B1)] 
$ H $ is a conservative Hunt process with values in $ S $; 
\item[(B2)] 
the state $ 0 $ is regular for itself, i.e., 
$ \bP_0(T_0=0)=1 $; 
\item[(B3)] 
the state $ 0 $ is recurrent, i.e., 
$ \bP_x(T_0<\infty)=1 $ for all $ x \in S $. 
\end{description}

Let $ H $ be an extension of $ H^0 $ satisfying Conditions {\bf (B1)}-{\bf (B3)}. 
Then the following assertions hold: 
\begin{enumerate}
\item 
there exists a positive continuous additive functional $ L=(L(t))_{t \ge 0} $ such that 
\\
$ \int_0^{\infty } 1_{\{ X(s) \neq 0 \}} \d L(s) = 0 $ 
(the process $ L $ is called the {\em local time of $ 0 $ for $ X $}); 
\item 
if $ A=(A(t))_{t \ge 0} $ is a non-negative continuous additive functional such that 
\\
$ \int_0^{\infty } 1_{\{ X(s) \neq 0 \}} \d A(s) = 0 $, 
then $ A(t) \equiv k L(t) $ for some constant $ k $. 
\end{enumerate}
For the proof of these facts, see, e.g., \cite[Theorem V.3.13]{BG}. 

We fix $ L $ for a choice of the local time of $ 0 $. 
We then see that there exists a constant $ \varsigma \ge 0 $ such that 
\begin{align}
\int_0^t 1_{\{ X(s) = 0 \}} \d s = \varsigma L(t) 
, \quad t \ge 0 . 
\label{eq: stag}
\end{align}
The constant $ \varsigma $ is called the {\em stagnancy rate}. 
Denote 
\begin{align}
\eta(l) = \inf \{ t \ge 0 : L(t) > l \} 
, \quad l \ge 0 . 
\label{}
\end{align}
For $ l \ge 0 $, we define $ p^{(l)} = (p^{(l)}(t))_{t \ge 0} \in D $ by 
\begin{align}
p^{(l)}(t) = 
\begin{cases}
X(\eta(l-)+t) 
& \text{if $ 0 \le t < \eta(l)-\eta(l-) $}, \\
0 
& \text{if $ t \ge \eta(l)-\eta(l-) $}. 
\end{cases}
\label{}
\end{align}
The point process $ p=(p^{(l)})_{l \ge 0} $ thus obtained 
will be called the {\em point process of excursions for $ \{ X,\bP_0 \} $}. 
It is then known (see \cite[Section 6]{MR0402949}) 
that $ \{ p,\bP_0 \} $ is a Poisson point process on $ D $ outside $ o $. 
Its characteristic measure will be denoted by $ \vn $ 
and called the {\em excursion measure}. 
We now see that 
\begin{align}
L = L(p,\varsigma) 
, \quad 
\eta = \eta(p,\varsigma) 
\quad \text{and} \quad 
X = X(p,\varsigma) . 
\label{}
\end{align}

\begin{Thm}[It\^o] \label{em}
Let $ (\vn,\varsigma) $ be as above. Then the following assertions hold: 
\begin{enumerate}
\item 
$ (\vn,\varsigma) $ satisfies Conditions {\bf (N0)}-{\bf (N3)}; 
\item 
for any $ t \ge 0 $, any $ A \in \cF_t $ and any $ A' \in \cB(D) $, it holds that 
\begin{align}
\vn(\{ T_0>t \} \cap A \cap \theta_t^{-1} A' ) 
= \vn \sbra{ \bP^0_{X(t)}(A') ; \{ T_0>t \} \cap A } , 
\label{eq: Mp}
\end{align}
provided that $ \vn(\{ T_0>t \} \cap A) < \infty $. 
\end{enumerate}
\end{Thm}

For the proof of Theorem \ref{em}, see It\^o \cite[Section 6]{MR0402949} 
and also Salisbury \cite{MR859837}. 

We also have the strong Markov property for $ \vn $ stated as follows. 

\begin{Thm} \label{sMp}
For any stopping time $ T $, 
any $ A \in \cF_T $ and any $ A' \in \cB(D) $, it holds that 
\begin{align}
\vn(\{ T_0>T \} \cap A \cap \theta_T^{-1} A' ) 
= \vn \sbra{ \bP^0_{X(T)}(A') ; \{ T_0>T \} \cap A } , 
\label{}
\end{align}
provided that $ \vn(\{ T_0>T \} \cap A) < \infty $. 
\end{Thm}

From this theorem we obtain the following corollary. 

\begin{Cor} \label{vnPo}
For any $ x \neq 0 $ and any $ A \in \cB(D) $, it holds that 
\begin{align}
\vn(\{ T_x<T_0 \} \cap \theta_{T_x}^{-1} A) 
= \vn(T_x<T_0) \bP^0_x(A) . 
\label{eq: vnbPo}
\end{align}
\end{Cor}

\Proof{
By Condition {\bf (N3)}, we have 
\begin{align}
\vn(\{ T_x<T_0 \}) \le \vn(\| X \| \ge | x |) < \infty . 
\label{eq: TxTo}
\end{align}
Hence we may apply Theorem \ref{sMp} for $ T=T_x $. 
Since $ X(T_x)=x $, we obtain \eqref{eq: vnbPo}. 
}

\subsection{Scaling property}

Let $ H $ be an extension of $ H^0 $ satisfying Conditions {\bf (B1)}-{\bf (B3)}. 
Let $ c>1 $ be a fixed constant. 
For $ \gamma>0 $, we define transformations $ \Psi_{\gamma} $ 
and $ \hat{\Psi}_{\gamma} $ of $ D $ by 
\begin{align}
(\Psi_{\gamma } w)(t) = c^{-\gamma } w(ct) 
, \quad 
(\hat{\Psi}_{\gamma } w)(t) = c^{-1} w(c^{\gamma}t) . 
\label{}
\end{align}
We introduce the following set of conditions: 
\begin{description}
\item[(S0)] 
$ c>1 $, $ \alpha >0 $, $ 0<\kappa <1/\alpha $ and $ c^{-\alpha } S \subset S $; 
\item[(S1)] 
$ \{ \Psi_{\alpha } X,\bP_x \} \law \{ X,\bP_{c^{-\alpha } x} \} $ 
for all $ x \in S $; 
\item[(S2)] 
$ \{ \Psi_{\alpha \kappa } L,\bP_0 \} \law \{ L,\bP_0 \} $. 
\end{description}

We need the following lemma. 

\begin{Lem} \label{stag0}
Suppose that Conditions {\bf (S0)}-{\bf (S2)} are satisfied. 
Then the stagnancy rate of the process $ \{ X,\bP_0 \} $ is necessarily equal to 0. 
\end{Lem}

The proof of Lemma \ref{stag0} will be given in Section \ref{sec: prsca}. 

Condition {\bf (S2)} is equivalent to the scaling property 
of the excursion measure as follows. 

\begin{Prop} \label{excscale}
Suppose that Conditions {\bf (S0)}-{\bf (S1)} are satisfied. 
Then Condition {\bf (S2)} is equivalent to the following condition: 
\begin{description}
\item[(S2)$ ' $] 
$ \vn \circ \Psi_{\alpha }^{-1} = c^{-\alpha \kappa } \vn $. 
\end{description}
\end{Prop}

\subsection{Homogenization theorem for jumping-in extensions}

In addition to Conditions {\bf (B1)}-{\bf (B3)}, 
we introduce the following set of conditions: 
\begin{description}
\item[(B4)] 
excursions leave $ 0 $ continuously, i.e., $ X(0)=0 $, $ \vn $-a.e.; 
\item[(B5)] 
excursions hit $ 0 $ continuously, i.e., $ X(T_0-)=0 $, $ \vn $-a.e. 
\end{description}

If Conditions {\bf (B1)}-{\bf (B5)} and {\bf (S0)}-{\bf (S2)} are satisfied, 
then we see, by Condition {\bf (S2)$ ' $}, that it also satisfies 
\begin{align}
\sigma(c^{-\alpha } x) = c^{\alpha \kappa } \sigma(x) 
\quad \text{for $ x \in S $}. 
\label{eq: psisigma}
\end{align}

Let $ H^0 $ be a Hunt process stopped upon hitting 0 
and let $ c>1 $, $ \alpha >0 $ and $ 0 < \kappa < 1/\alpha $ be fixed. 
Let $ S' $ be a measurable space 
and let $ \{ \vn_v \}_{v \in S'} $ be a kernel on $ D $. 
We introduce the following condition: 
\begin{description}
\item[(B)] 
for each $ v \in S' $, 
the measure $ \vn_v $ 
is the excursion measure of an extension $ H_v = \{ X,(\bP^v_x)_{x \in S} \} $ 
of $ H^0 $ satisfying Conditions {\bf (B1)}-{\bf (B5)} and {\bf (S0)}-{\bf (S2)}. 
\end{description}

For a finite measure $ \rho $ on $ S' $ 
and a $ \sigma $-finite measure $ j $ on $ S \setminus \{ 0 \} $, we define 
\begin{align}
\vn_{\rho,j}(\d w) = \int_{S'} \rho(\d v) \vn_v(\d w) 
+ \int_{S \setminus \{ 0 \}} j(\d x) \bP^0_x(\d w) . 
\label{eq: vnthetaj}
\end{align}
For a triplet $ (\rho,j,\varsigma) $, 
we introduce the following condition: 
\begin{description}
\item[(C1)] 
the pair $ (\vn_{\rho,j},\varsigma) $ satisfies Conditions {\bf (N0)}-{\bf (N3)}. 
\item[(C2)] 
there exists a measurable map $ \psi:S \setminus \{ 0 \} \to S' $ such that, 
for $ j $-a.e. $ x \in S \setminus \{ 0 \} $, 
$ \sigma_{\psi(x)}(x)>0 $ 
and 
$ \psi(c^{-\alpha n}x) = \psi(x) $ for all $ n \in \bN \cup \{ \infty \} $. 
\end{description}
For a triplet $ (\rho,j,\varsigma) $ satisfying Condition {\bf (C)}, 
let $ \{ p_{\rho,j},\bP \} $ be a Poisson point process on $ D $ outside $ o $ 
with characteristic measure $ \vn_{\rho,j} $. 
We write 
\begin{align}
X_{\rho,j,\varsigma}(t) = X(p_{\rho,j},\varsigma;t) , 
\quad 
L_{\rho,j,\varsigma}(t) = L(p_{\rho,j},\varsigma;t) , 
\quad 
\eta_{\rho,j,\varsigma}(l) = \eta(p_{\rho,j},\varsigma;l) 
\label{}
\end{align}
and call $ \{ X_{\rho,j,\varsigma},\bP \} $ a {\em jumping-in extension} 
of the minimal process $ H^0 $.

For a scaling exponent $ \gamma>0 $, we define 
\begin{align}
\vn^{(n)}_{\rho,j} 
= c^{\gamma n} \vn_{\rho,j} \circ (\Psi_{\alpha }^n)^{-1} 
, \quad 
\varsigma^{(n)} 
= c^{- (1-\gamma) n} \varsigma 
\label{eq: scaled1}
\end{align}
and 
\begin{align}
X_{\rho,j,\varsigma}^{(n)} = 
\Psi_{\alpha }^n X_{\rho,j,\varsigma} , 
\quad 
L_{\rho,j,\varsigma}^{(n)} = 
\Psi_{\gamma}^n L_{\rho,j,\varsigma} , 
\quad 
\eta_{\rho,j,\varsigma}^{(n)} = 
\hat{\Psi}_{\gamma }^n \eta_{\rho,j,\varsigma} . 
\label{eq: scaled5}
\end{align}

Let $ H^0 $ be a Hunt process stopped upon hitting 0 
and let $ c>1 $, $ \alpha>0 $ and $ 0<\kappa < 1/\alpha $ be constants. 
Let $ \{ \vn_v \}_{v \in S'} $ be a kernel satisfying Condition {\bf (B)}. 
Denote 
\begin{align}
\sigma_v(x) = \vn_v(T_x<T_0) \quad \text{for $ v \in S' $}. 
\label{}
\end{align}
Let $ (\rho,j,\varsigma ) $ satisfy Conditions {\bf (C1)}-{\bf (C2)}. 
In order to handle various examples together, 
we give the following two auxiliary theorems. 

\begin{Thm}[jumping-in vanishing case] \label{hom}
Suppose the following condition: 
\begin{description}
\item[(C3)] 
$ T_{c^{-\alpha n} x} \to 0 $, $ \vn_{\psi(x)} $-a.e. 
for $ j $-a.e. $ x \in S \setminus \{ 0 \} $; 
\item[(C4)] 
$ (\vn_{\rho^*,0},0) $ satisfies Conditions {\bf (N0)}-{\bf (N3)}, where 
$ \rho^* $ is the finite measure on $ S' $ defined by 
\begin{align}
\rho^* = \rho + 
\int_{S \setminus \{ 0 \}} \frac{j(\d x)}{\sigma_{\psi(x)}(x)} \delta_{\psi(x)} , 
\label{}
\end{align}
\end{description}
where $ \delta $ deonte the Dirac delta. 
Then, for the scaling exponent $ \gamma = \alpha \kappa $, 
it holds as $ n \to \infty $ that 
\begin{align}
\vn^{(n)}_{\rho,j} \to \vn_{\rho^*,0} , 
\quad 
\varsigma^{(n)} \to 0 
\label{}
\end{align}
and 
\begin{align}
\cbra{ \rbra{ X_{\rho,j,\varsigma}^{(n)} , L_{\rho,j,\varsigma}^{(n)} , 
\eta_{\rho,j,\varsigma}^{(n)} } , \bP } 
\claw \cbra{ \rbra{ X_{\rho^*,0,0} , L_{\rho^*,0,0} , \eta_{\rho^*,0,0} \Big. } , \bP } . 
\label{}
\end{align}
\end{Thm}

\begin{Thm}[jumping-in dominant case] \label{hom2}
Suppose the following condition: 
\begin{description}
\item[(C5)] 
there exist a $ \sigma $-finite measure $ \mu $ on a measurable space $ S'' $, 
measurable mappings $ J_n:S'' \to S \setminus \{ 0 \} $, 
$ J^*:S'' \to S \setminus \{ 0 \} $, 
and a constant $ \beta \in (0,\kappa) $ 
such that 
\begin{align}
c^{\alpha \beta n} 
\int_{S \setminus \{ 0 \}} 
\frac{j(\d x)}{\sigma_{\psi(x)}(c^{-\alpha n} x)} f(\psi(x),c^{-\alpha n} x) 
= \int_{S''} \mu(\d y) f(\psi(J^*y),J_ny) 
\label{}
\end{align}
for all $ n \in \bN $ 
and all non-negative measurable function $ f $ on $ S' \times S $ 
and 
\begin{align}
T_{J_n y} \to T_{J^* y} 
\ \text{and} \ 
X(T_{J^*y}-) = X(T_{J^*y}) 
, \ \text{$ \vn_{\psi(J^* y)} $-a.e. for $ \mu $-a.e. $ y \in S'' $}; 
\label{}
\end{align}
\item[(C6)] 
$ (\vn_{0,j^*},0) $ satisfies Conditions {\bf (N0)}-{\bf (N3)}, where 
$ j^* $ is the $ \sigma $-finite measure on $ S \setminus \{ 0 \} $ defined by 
\begin{align}
j^* = \int_{S''} \mu(\d y) \sigma_{\psi(J^*y)}(J^*y) \delta_{J^*y} . 
\label{}
\end{align}
\end{description}
Then, for the scaling exponent $ \gamma = \alpha \beta $, 
it holds as $ n \to \infty $ that 
\begin{align}
\vn^{(n)}_{\rho,j} \to \vn_{0,j^*} 
, \quad 
\varsigma^{(n)} \to 0 
\label{}
\end{align}
and 
\begin{align}
\cbra{ \rbra{ X_{\rho,j,\varsigma}^{(n)} , L_{\rho,j,\varsigma}^{(n)} , 
\eta_{\rho,j,\varsigma}^{(n)} } , \bP } 
\claw \cbra{ \rbra{ X_{0,j^*,0} , L_{0,j^*,0} , \eta_{0,j^*,0} \Big. } , \bP } . 
\label{}
\end{align}
\end{Thm}

\section{Examples} \label{sec: ex}

\subsection{Positive self-similar Markov processes}

Let $ \alpha $ and $ \kappa $ be positive numbers 
such that $ 0<\kappa <1/\alpha $ 
and let $ c>1 $ be an arbitrary number. 
Let $ \{ X,(\bP_x)_{x \ge 0} \} $ be a Hunt process with values in $ S=[0,\infty ) $ such that 
{\bf (B1)}-{\bf (B5)} and {\bf (S0)}-{\bf (S2)} hold and 
\begin{align}
\begin{cases}
T_{x_n} \to 0, \ \text{$ \vn $-a.e.} & \text{if $ x_n \to 0 $}, \\
T_{x_n} \to T_x, \ \text{$ \vn $-a.e.} & \text{if $ x_n \to x > 0 $}, 
\end{cases}
\label{eq: pssMp Tx conv}
\end{align}
where $ \vn $ denotes the excursion measure away from 0 
according to a particular choice of the local time at 0. 
Let $ \bP^0_x $ denote the law of $ X(t \wedge T_0) $ under $ \bP_x $. 
Then Condition {\bf (B)} is satisfied 
for $ S' = \{ 0 \} $, $ H^0 = \{ X,(\bP^0_x)_{x>0} \} $ and $ \vn_0 = \vn $. 

Such a process can be obtained in the following manner. 
Let $ \{ Z,(\bQ_z)_{z \in \bR} \} $ be a L\'evy process which satisfies 
the following conditions: 
\begin{description}
\item[(P1)] 
$ Z $ drifts to $ -\infty $; 
\item[(P2)] 
$ Z $ is spectrally negative; 
\item[(P3)] 
every point is regular for itself; 
\item[(P4)] 
every point is accessible. 
\end{description}
By {\bf (P1)} and {\bf (P2)}, it is known (see, e.g., \cite[Section 8.1]{Kyp}) that 
there exists a unique constant $ \kappa>0 $ such that 
the following {\em Cram\'er condition} is satisfied: 
\begin{align}
\bQ_0 \sbra{ \exp(\kappa Z(u)) } = 1 
\quad \text{for all $ u \ge 0 $}. 
\label{eq: cramer}
\end{align}
We recall the {\em Lamperti transformaiton} following Lamperti \cite{MR0307358} as follows. 
Let $ \alpha $ be a fixed constant such that $ 0 < \alpha < 1/\kappa $. 
Define 
\begin{align}
\tau(u) = \int_0^u \exp(Z(s)/\alpha ) \d s 
\label{}
\end{align}
and 
\begin{align}
Y(t) = 
\begin{cases}
\exp(Z(\tau^{-1}(t)) & \text{for $ 0 \le t < \tau(\infty ) $}, \\
0 & \text{for $ t \ge \tau(\infty ) $}. 
\end{cases}
\label{}
\end{align}
For $ x>0 $, we write $ \bP^0_x $ the law of $ Y $ under $ \bQ_{\log x} $ 
and let $ H^0 = \{ X,(\bP^0_x)_{x>0} \} $. 
By the theorem obtained by Rivero \cite{MR2146891,MR2364226}
and Fitzsimmons \cite{MR2266714} independently, 
we see, thanks to the Cram\'er condition \eqref{eq: cramer}, that 
there exists a unique $ \alpha $-self-similar recurrent extension of $ H^0 $ 
whose excursions leave $ 0 $ continuously, 
which we will denote by $ \{ X,(\bP_x)_{x \ge 0} \} $. 
Then we see that {\bf (B1)}-{\bf (B5)} and {\bf (S0)}-{\bf (S2)} are satisfied. 
Since $ z \mapsto T_{z+}(Z) $ is a subordinator and has no fixed discontinuity, 
we see that $ \bP^0_{\eps}(D \setminus \{ T_{x_n} \to T_x \}) = 0 $ 
for any sequence $ \{ x_n \} $ converging to $ x>\eps>0 $. 
Using Corollary \ref{vnPo}, we see that $ T_{x_n} \to T_x $, $ \vn $-a.e. 
for any sequence $ \{ x_n \} $ converging to $ x>0 $. 
Since $ \{ X,\vn \} $ has c\`adl\`ag paths and has no positive jumps, 
we further see that $ T_{x_n} \to 0 $, $ \vn $-a.e. 
for any sequence $ \{ x_n \} $ converging to 0. 
Consequently we have verified that \eqref{eq: pssMp Tx conv} is satisfied. 
We have thus obtained $ \{ X,(\bP_x)_{x \ge 0} \} $ as desired. 

Since \eqref{eq: psisigma} holds for any $ c>0 $, we see that 
\begin{align}
\sigma(x) := \vn(T_x < T_0) = \delta x^{-\kappa } 
\label{eq: pssMp sigma}
\end{align}
for $ \delta = \vn(T_1<T_0) $. 
We need the following. 

\begin{Lem} \label{pssMp PT0w1} 
There exist positive constants $ c_1 $ and $ c_2 $ such that 
\begin{align}
c_1 (x^{\kappa} \wedge 1) \le \bP^0_x [T_0 \wedge 1] \le c_2 (x^{\kappa} \wedge 1) 
\quad \text{for all $ x>0 $}. 
\label{}
\end{align}
\end{Lem}

\Proof{
For $ x \ge 1 $, we see by the scaling property that 
\begin{align}
0 < \bP^0_1(T_0 \ge 1) 
\le \bP^0_1(T_0 \ge x^{-1/\alpha }) 
= \bP^0_x(T_0 \ge 1) \le \bP^0_x[T_0 \wedge 1] \le 1 . 
\label{}
\end{align}
By Corollary \ref{vnPo} and by \eqref{eq: pssMp sigma}, we have 
\begin{align}
\bP^0_x[T_0 \wedge 1] = \delta^{-1} x^{\kappa } 
\vn \sbra{ T_0 \circ \theta_{T_x} \wedge 1 ; T_x < T_0 } . 
\label{}
\end{align}
For $ 0<x<1 $, 
we see by {\bf (P2)} that $ \{ T_1<T_0 \} \subset \{ T_x < T_1 < T_0 \} $, $ \vn $-a.e., 
and hence we obtain 
\begin{align}
0 < \vn \sbra{ T_0 \circ \theta_{T_1} \wedge 1 ; T_1 < T_0 } 
\le \vn \sbra{ T_0 \circ \theta_{T_x} \wedge 1 ; T_x < T_0 } 
\le \vn [T_0 \wedge 1] < \infty . 
\label{}
\end{align}
The proof is now complete. 
}

We identify a measure $ \rho $ on $ S'=\{ 0 \} $ with a positive number $ \rho(\{ 0 \}) $. 
By Lemma \ref{pssMp PT0w1}, we see that, 
a pair $ (\vn_{\rho,j},\varsigma) $ satisfies {\bf (N0)}-{\bf (N3)} 
if and only if the following conditions are satisfied: 
\begin{description}
\item[(CP1)] 
$ \rho $ is a non-negative constant; 
\item[(CP2)] 
$ j $ satisfies $ \textstyle \int_{(0,\infty )} (x^{\kappa } \wedge 1) j(\d x) < \infty $; 
\item[(CP3)] 
any one of the following holds: 
$ \rho>0 $, $ j((0,\infty )) = \infty $ and $ \varsigma>0 $. 
\end{description}

\begin{Cor}[jumping-in vanishing case] \label{psj}
Let $ (\rho,j,\varsigma) $ satisfy {\bf (CP1)}-{\bf (CP3)}. 
Suppose, moreover, that 
\begin{align}
\rho^* := \rho + \frac{1}{\delta} \int_{(0,\infty )} x^{\kappa} j(\d x) < \infty . 
\label{}
\end{align}
Then the same assertions as Theorem \ref{hom} hold. 
\end{Cor}

Corollary \ref{psj} is an immediate consequence of Theorem \ref{hom}, 
and so we omit its proof. 

\begin{Cor}[jumping-in dominant case] \label{pbj}
Let $ (\rho,j,\varsigma) $ satisfy {\bf (CP1)}-{\bf (CP3)}. 
Let $ \beta \in (0,\kappa ) $ and $ j_0>0 $ be constants. 
Suppose, moreover, that 
\begin{align}
j((x,\infty )) \sim j_0 x^{-\beta} 
\quad \text{as $ x \to \infty $}. 
\label{eq: regvar j}
\end{align}
Define a $ \sigma $-finite measure $ j^* $ 
on $ (0,\infty ) $ by 
\begin{align}
j^*(\d x) = j_0 \beta x^{-\beta-1} \d x . 
\label{}
\end{align}
Then the same assertions as Theorem \ref{hom2} hold. 
\end{Cor}

\Proof[Proof of Corollary \ref{pbj}]{
We define a function $ J:(0,\infty ) \to (0,\infty ) $ by 
\begin{align}
J(y) = \inf \cbra{ x>0 : \frac{1}{\delta} \int_{(0,x]} s^{\kappa} j(\d s) > y } . 
\label{eq: J}
\end{align}
In the same way as $ J $ we define $ J^* $ with $ j $ being replaced by $ j^* $. 
Set $ S''=(0,\infty ) $, $ \mu(\d y) = \d y $ 
and $ J_ny = c^{-\alpha n} J(c^{\alpha (\kappa -\beta)n}y) $. 
By \eqref{eq: regvar j}, 
we see that 
\begin{align}
\frac{1}{\delta} \int_{(0,x]} s^{\kappa } j(\d s) 
\sim \frac{j_0 \beta}{\delta (\kappa -\beta)} x^{\kappa - \beta} 
= \frac{1}{\delta} \int_{(0,x]} s^{\kappa } j^*(\d s) 
\quad \text{as $ x \to \infty $}, 
\label{}
\end{align}
which shows that {\bf (C5)} and {\bf (C6)} are satisfied. 
We can thus apply Theorem \ref{hom2}. 
}

\subsection{Walsh diffusions}

Let us take up the {\em Walsh diffusions}, 
which have been first introduced in Walsh \cite[Epilogue]{WalshBM} 
and developed in Barlow--Pitman--Yor \cite{MR1022917}. 
A Walsh diffusion is a diffusion process with values in $ \bR^2 $ 
whose stopped process starting from $ x \neq 0 $ and stopped at $ 0 $ 
takes values in the ray $ R(x) := \{ rx: r \ge 0 \} $. 
In this paper we confine ourselves to the case where 
the stopped processes are Bessel ones. 

Let $ 0 < \alpha < 1 $ and $ c>1 $ be fixed. 
Let $ B = \{ X,(\bQ_r)_{r \ge 0} \} $ 
denote the $ \frac{\d}{\d m} \frac{\d}{\d x} $-diffusion process 
for $ m(x) = x^{1/\alpha -1} $ 
which takes values in $ [0,\infty ) $ 
and which has 0 as the reflecting boundary. 
Let $ \vn_B $ denote the excursion measure of $ B $ away from 0 according 
to a particular choice of the local time at 0 which we will denote by $ L_B $. 
Then {\bf (S0)}-{\bf (S2)} holds for $ S=[0,\infty ) $ and $ \kappa = 1 $ 
(see, e.g., \cite{MR2417969}). 
Let $ \bQ^0_r $ denote the law of $ X(t \wedge T_0) $ under $ \bQ_r $. 
By the same argument as the proof of Lemma \ref{pssMp PT0w1} 
(see also \cite[Example 6.1]{ItoKyoto}), we see that 
\begin{align}
c_1 (r \wedge 1) \le \bQ_r[T_0 \wedge 1] \le c_2 (r \wedge 1) 
\quad \text{for all $ r>0 $} 
\label{eq: QrT0}
\end{align}
holds for some positive constants $ c_1 $ and $ c_2 $. 

Let $ S=\bR^2 $ and let $ S'=S^1 = \{ v \in \bR^2 : |v|=1 \} $. 
Define $ \psi:\bR^2 \setminus \{ 0 \} \to S^1 $ by $ \psi(x) = x/|x| $. 
We define $ (\vn_v)_{v \in S^1} $ 
and $ (\bP^0_x)_{x \in \bR^2 \setminus \{ 0 \}} $ by 
\begin{align}
\vn_v = \int_{D_{[0,\infty )}} \vn_B(\d q) \delta_{qv} 
\quad \text{and} \quad 
\bP^0_x = \int_{D_{[0,\infty )}} \bQ^0_{|x|}(\d q) \delta_{q\psi(x)} , 
\label{}
\end{align}
where for $ q = (q(t))_{t \ge 0} \in D_{[0,\infty )} $ 
we write $ qv = (q(t)v)_{t \ge 0} \in D_{\bR^2} $. 
Since $ \vn_B(T_r<T_0) = \delta r^{-1} $ for all $ r>0 $ 
with $ \delta = \vn_B(T_1<T_0) $, we see that 
\begin{align}
\vn_{\psi(x)}(T_x<T_0) = \delta |x|^{-1} 
, \quad x \in \bR^2 \setminus \{ 0 \} . 
\label{}
\end{align}
For $ v \in S^1 $ and for $ \{ x_n \} \subset R(v) $, 
we easily see that 
\begin{align}
\begin{cases}
T_{x_n} \to 0, \ \text{$ \vn_v $-a.e.} & \text{if $ x_n \to 0 $}, \\
T_{x_n} \to T_x, \ \text{$ \vn_v $-a.e.} & \text{if $ x_n \to x \neq 0 $}. 
\end{cases}
\label{eq: Wd Tx conv}
\end{align}
For $ v \in S^1 $ and $ x \in R(v) $, we define 
$ \bP^v_x = \int_{D_{[0,\infty )}} \bQ_{|x|}(\d q) \delta_{qv} $. 
For $ v \in S^1 $ and $ x \notin R(v) $, 
we define $ \bP^v_x $ as the law of the process 
which is obtained as $ X(t \wedge T_0(X)) + \tilde{X}((t-T_0(X)) \vee 0) $, 
where $ \{ X,\bP^0_x \} $ and $ \{ \tilde{X},\bP_0 \} $ are independent processes 
defined on a common probability space. 
Then $ H_v = \{ X,(\bP^v_x)_{x \in \bR^2} \} $ 
is an extension of $ H^0 = \{ X,(\bP^0_x)_{x \in \bR^2} \} $ 
and $ \vn_v $ is the excursion measure away from 0 of $ H_v $. 
It is immediate that {\bf (B)} holds for $ \kappa =1 $.

For measures $ \rho $ on $ S^1 $ and $ j $ on $ \bR^2 \setminus \{ 0 \} $, 
we define $ \vn_{\rho,j} $ on $ D_{\bR^2} $ by \eqref{eq: vnrhoj}. 
By \eqref{eq: QrT0}, we see that 
a pair $ (\vn_{\rho,j},\varsigma) $ satisfies {\bf (N0)}-{\bf (N3)} 
if and only if the following conditions are satisfied: 
\begin{description}
\item[(CW1)] 
$ \rho $ is a finite measure on $ S^1 $; 
\item[(CW2)] 
$ j $ satisfies $ \int_{\bR^2 \setminus \{ 0 \}} (|x| \wedge 1) j(\d x) < \infty $; 
\item[(CW3)] 
any one of the following holds: 
$ \rho(S^1)>0 $, $ j(\bR^2 \setminus \{ 0 \}) = \infty $ and $ \varsigma>0 $. 
\end{description}
We identify $ S^1 \times (0,\infty ) $ with $ \bR^2 \setminus \{ 0 \} $ 
via the bicontinuous bijection $ (v,r) \mapsto rv $. 
A $ \sigma $-finite measure $ j(\d v \d r) $ on $ S^1 \times (0,\infty ) $ 
allows at least one disintegration of the form 
\begin{align}
j(\d v \d r) = \rho_j(\d v) j_v(\d r) 
\label{eq: disinteg}
\end{align}
for a finite measure $ \rho_j $ on $ S^1 $ and 
a kernel $ \{ j_v \}_{v \in S^1} $ on $ (0,\infty ) $. 
We can obtain such a disintegration, for example, as follows: 
Take a measurable function $ f(v,r) $ which is positive $ j $-a.e. 
and which satisfies $ j[f]=1 $. 
Then, by conditioning, the probability measure $ \bar{j}(\d v \d r) = f(v,r) j(\d v \d r) $ 
possesses a unique disintegration 
$ \bar{j} = \rho_{\bar{j}}(\d v) \bar{j}_v(\d r) $ 
with a probability measure $ \rho_{\bar{j}} $ 
and a probability kernel $ \{ \bar{j}_v \}_{v \in S^1} $. 
We set $ \rho_j = \rho_{\bar{j}} $ 
and $ j_v(\d r) = f(v,r)^{-1} \bar{j}_v(\d r) $ 
and then we obtain the disintegration \eqref{eq: disinteg}. 
The disintegration \eqref{eq: disinteg} is not unique; 
in fact, for any bounded measurable function $ \pi $ on $ S^1 $ with positive values, 
we obtain another disintegration $ j(\d v \d r) = \rho'_j(\d v) j'_v(\d r) $, 
where $ \rho'_j(\d v) = \pi(v) \rho_j(\d v) $ 
and $ j'_v(\d r) = \pi(v)^{-1} j_v(\d r) $.

\begin{Cor}[jumping-in vanishing case] \label{Walsh1}
Let $ (\rho,j,\varsigma) $ satisfy {\bf (CW1)}-{\bf (CW3)}. 
Let $ j(\d v \d r) = \rho_j(\d v) j_v(\d r) $ be a disintegration of $ j $ 
and suppose that 
\begin{align}
\int_{\bR^2 \setminus \{ 0 \}} |x| j(\d x) 
= \int_{S^1} \pi(v) \rho_j(\d v) < \infty , 
\label{}
\end{align}
where $ \pi(v) = \int_{(0,\infty )} r j_v(\d r) $. 
Define the finite measure $ \rho^* $ on $ S^1 $ by 
\begin{align}
\rho^*(\d v) = 
\rho(\d v) + \frac{1}{\delta} \pi(v) \rho_j(\d v) . 
\label{}
\end{align}
Then the same assertions as Theorem \ref{hom} hold. 
\end{Cor}

Corollary \ref{Walsh1} is an immediate consequence of Theorem \ref{hom}, 
and so we omit its proof.

\begin{Cor}[jumping-in dominant case] \label{Walsh2}
Let $ (\rho,j,\varsigma) $ satisfy {\bf (CW1)}-{\bf (CW3)}. 
Let $ j(\d v \d r) = \rho_j(\d v) j_v(\d r) $ be a disintegration of $ j $ 
and suppose that 
there exist a constant $ 0 < \beta < 1 $ 
and a non-negative measurable function $ \pi $ on $ S^1 $ 
such that 
\begin{align}
\int_{S^1} \pi(v) \rho_j(\d v) \in (0,\infty ) 
\label{}
\end{align}
and, for any $ v \in S^1 $, 
\begin{align}
r^{\beta} j_v((r,\infty )) \to \pi(v) 
\quad \text{as $ r \to \infty $}. 
\label{eq: regvar jv}
\end{align}
Define a $ \sigma $-finite measure $ j^* $ 
on $ \bR^2 \setminus \{ 0 \} \simeq S^1 \times (0,\infty ) $ by 
\begin{align}
j^*(\d v \d r) = \pi(v) \rho_j(\d v) \beta r^{-\beta-1} \d r . 
\label{}
\end{align}
Then the same assertions as Theorem \ref{hom2} hold. 
\end{Cor}

\Proof[Proof of Corollary \ref{Walsh2}]{
Note that $ j^* $ admits a disintegration $ j^*(\d v \d r) = \rho_{j^*}(\d v) j^*_v(\d r) $ 
where $ \rho_{j^*}(\d v) = \pi(v) \rho_j(\d v) $ 
and $ j^*_v(\d r) = \beta r^{-\beta-1} \d r $. 

Let $ S'' = S^1 \times (0,\infty ) $ 
and let $ \mu(\d v \d r) = \d r $. 
For $ v \in S^1 $, we define 
a function $ J^0_v:(0,\infty ) \to (0,\infty ) $ by 
\begin{align}
J^0_vy = \inf \cbra{ r>0 : \frac{1}{\delta} \int_{(0,r]} s j_v(\d s) > y } 
\label{}
\end{align}
and define a function $ J:S^1 \times (0,\infty ) \to \bR^2 \setminus \{ 0 \} $ 
by $ J(v,y) = (J^0_vy)v $. 
In the same way as $ J $ we define $ J^* $ 
with $ j_v $ being replaced by $ j^*_v $. 
For $ n \in \bN $, set 
$ J_n(v,y) = c^{-\alpha n} J(v,c^{\alpha (1-\beta) n}y) $. 
By the assumption \eqref{eq: regvar jv}, we have 
\begin{align}
\frac{1}{\delta} \int_{(0,r]} s j_v(\d s) 
\sim \frac{\pi(v) \beta}{\delta(1-\beta)} r^{1-\beta} 
= \frac{1}{\delta} \int_{(0,r]} s j^*_v(\d s) 
\quad \text{as $ r \to \infty $}, 
\label{}
\end{align}
which shows that {\bf (C5)} and {\bf (C6)} are satisfied. 
We can thus apply Theorem \ref{hom2}. 
}

\subsection{The Brownian motion on the Sierpi\'nski gasket}

We take up the Brownian motion on the {\em Sierpi\'nski gasket}. 
For its precise definition and several facts which we will utilize later, 
see Barlow--Perkins \cite{MR966175}. 

Let $ S=G $ denote the Sierpi\'nski gasket in $ \bR^2 $ 
and let $ \{ X,(\bP_x)_{x \in G} \} $ denote the Brownian motion on $ G $. 
Noting that every point of $ G $ is regualr for itself, 
we let $ L_t^x $ denote a jointly continuous version of the local time 
and denote $ L(t) = L_t^0 $. 
Then {\bf (S0)}-{\bf (S2)} hold for 
\begin{align}
c = 5 
, \quad 
\alpha = \frac{\log 2}{\log 5} 
, \quad 
\kappa = \frac{\log 5 - \log 3}{\log 2} . 
\label{}
\end{align}
Let $ \vn $ denote the excursion measure away from 0. 

Let $ \bP^0_x $ denote the law of $ X(t \wedge T_0) $ under $ \bP_x $ 
and let $ H^0 = \{ X,(\bP^0_x)_{x \in G} \} $. 
We denote $ x = (x_1,x_2) \in \bR^2 $ 
and set $ G_{\pm} = \{ x \in G : \pm x_1 \ge 0 \} $. 
We write $ \bP^{\pm}_x $ for the law of $ Y_{\pm} $ under $ \bP_x $ 
and write $ H_{\pm} = \{ X,(\bP^{\pm}_x)_{x \in G} \} $, 
where 
\begin{align}
Y_{\pm}(t) = 
\begin{cases}
X(t) & \text{for $ 0 \le t \le T_0 $}, \\
(\pm |X_1(t)|,X_2(t)) & \text{for $ t > T_0 $}. 
\end{cases}
\label{}
\end{align}
We then see that $ H_{\pm} $ are extensions of $ H^0 $ 
whose excursion measures away from 0 are 
\begin{align}
\vn_{\pm} = \vn|_{\{ X(t) \in G_{\pm} \ \text{for all $ t \ge 0 $} \}} . 
\label{}
\end{align}
Note that $ \vn = \vn_+ + \vn_- $. 
Letting $ S' = \{ +,- \} $, we see that Condition {\bf (B)} is satisfied. 
We identify a measure $ \rho $ on $ S' $ 
with the pair $ (\rho_+,\rho_-) := (\rho(\{ + \}),\rho(\{ - \})) $. 
A pair $ (\vn_{\rho,j},\varsigma) $ satisfies {\bf (N0)}-{\bf (N3)} 
if and only if the following conditions are satisfied: 
\begin{description}
\item[(CG1)] 
$ \rho_+ $ and $ \rho_- $ are non-negative finite constants; 
\item[(CG2)] 
$ j $ satisfies $ \textstyle \int_{G \setminus \{ 0 \}} j(\d x) \bP^0_x[T_0 \wedge 1] < \infty $; 
\item[(CG4)] 
any one of the following holds: 
$ \rho_+>0 $, $ \rho_->0 $, $ j(G \setminus \{ 0 \}) = \infty $ and $ \varsigma>0 $. 
\end{description}

To obtain the homogenization theorem, we need the following. 

\begin{Lem} \label{kajino}
For $ \{ x_n \} \subset G $ such that $ x_n \to x \neq 0 $, 
it holds that 
$ T_{x_n} \to T_x $, $ \vn $-a.e. 
\end{Lem}

\Proof{
Suppose that the following assertion is established: 
\begin{align}
T_{x_n} \to T_x , \quad \text{$ \bP_a $-a.e. for all $ a \neq 0 $}. 
\label{eq: SG Txconv}
\end{align}
If $ T_x < T_0 $, then $ L_{T_0}^y > 0 $ for all $ y $ in some neighborhood of $ x $, 
so that we have $ T_{x_n} \wedge T_0 = T_{x_n} \to T_x = T_x \wedge T_0 $ as $ n \to \infty $. 
If $ T_x > T_0 $, then $ L_{T_0}^y = 0 $ for all $ y $ in some neighborhood of $ x $, 
so that we have $ T_{x_n} \wedge T_0 = T_x \wedge T_0 = T_0 $ for large $ n $. 
We thus see that $ \bP^0_a(D \setminus \{ T_{x_n} \to T_x \}) = 0 $ for all $ a \neq 0 $. 
By the Markov property of $ \vn $, we obtain the desired result. 

Let us now prove \eqref{eq: SG Txconv}\footnote{This proof is due to N. Kajino.}. 
\\ $ 1^{\circ}) $. 
Let us prove $ T_x \le \liminf T_{x_n} $. 
Suppose $ t_0 := \liminf T_{x_n} < \infty $. 
Then there exists a subsequence $ \{ n(k) \} $ such that $ T_{n(k)} \to t_0 $, 
so that we have $ X(t_0) = \lim X(T_{x_{n(k)}}) = \lim x_{n(k)} = x $. 
This shows $ T_x \le t_0 = \liminf T_{x_n} $. 
\\ $ 2^{\circ}) $. 
Let us prove $ T_x \ge \limsup T_{x_n} $. 
Suppose $ T_x < T_0 $. 
For $ T_x < t_0 < T_0 $, 
we have $ L^y_{t_0} > 0 $ for all $ y $ in some neighborhood of $ x $, 
so that we have $ \limsup T_{x_n} \le t_0 $. 
This shows $ \limsup T_{x_n} \le T_x $. 

The proof is now complete. 
}

For a $ \sigma $-finite measure $ j $ on $ G \setminus \{ 0 \} $, 
set $ j_+ = j|_{G_+ \setminus \{ 0 \}} $ and 
$ j_- = \check{j}|_{G_+ \setminus \{ 0 \}} $, 
where $ \check{j} $ is the pullback of $ j $ under $ (x_1,x_2) \mapsto (-x_1,x_2) $. 
We define mappings 
$ \phi_1 : G_+ \setminus \{ 0 \} \to [0,1] $ and 
$ \phi_2 : G_+ \setminus \{ 0 \} \to (0,\infty ) $ by 
\begin{align}
\phi_1(x_1,x_2) = \frac{2x_2}{\sqrt{3} x_1 + x_2} 
, \quad 
\phi_2(x_1,x_2) = x_1 + \frac{1}{\sqrt{3}} x_2 
, \quad (x_1,x_2) \in G_+ . 
\label{}
\end{align}
We then see that the mapping 
$ \phi = (\phi_1,\phi_2):G_+ \setminus \{ 0 \} \to [0,1] \times (0,\infty ) $ 
is a measurable injection. 
For $ v \in [0,1] $, we write 
\begin{align}
R(v) = \{ \phi_2(x) : x \in G_+ \setminus \{ 0 \} \ \text{and} \ \phi_1(x)=v \} . 
\label{}
\end{align}
We then see that $ c^{-\alpha } R(v) = R(v) $. 
Since the pullbacks $ j_{\pm} \circ (\phi_1,\phi_2)^{-1} $ 
are $ \sigma $-finite measures on $ [0,1] \times (0,\infty ) $, 
we may obtain at least one disintegration of the form 
\begin{align}
(j_{\pm} \circ (\phi_1,\phi_2)^{-1})(\d v \d r) 
= \rho^{\pm}_{j}(\d v) j^{\pm}_v(\d r) 
\label{}
\end{align}
for finite measures $ \rho^{\pm}_j $ on $ [0,1] $ and 
kernels $ \{ j^{\pm}_v \}_{v \in [0.1]} $ on $ (0,\infty ) $ such that 
$ j^{\pm}_v((0,\infty ) \setminus R(v)) = 0 $ for all $ v \in [0,1] $.

\begin{Cor}[jumping-in dominant case] \label{SG2}
Let $ (\rho,j,\varsigma) $ satisfy {\bf (CW1)}-{\bf (CW3)}. 
Let $ (j_{\pm} \circ (\phi_1,\phi_2)^{-1})(\d v \d r) 
= \rho^{\pm}_j(\d v) j^{\pm}_v(\d r) $ be disintegrations 
and suppose that 
there exist a constant $ 0 < \beta < 1 $ 
and a $ \sigma $-finite measure $ j^* $ on $ G \setminus \{ 0 \} $ 
with disintegrations $ (j^*_{\pm} \circ (\phi_1,\phi_2)^{-1})(\d v \d r) 
= \rho^{\pm}_{j^*}(\d v) j^{*,\pm}_v(\d r) $ 
such that 
$ j^{*,\pm}_v((0,\infty ) \setminus R(v)) = 0 $ for any $ v \in [0,1] $ 
and, for any $ v \in [0,1] $ and any $ r \in R(v) $, 
\begin{align}
c^{- \alpha (\kappa -\beta ) n} 
\int_{(0,c^{\alpha n} r]} \frac{j^{\pm}_v(\d s)}{\sigma(\phi^{-1}(s,v))} 
\to \int_{(0,r]} \frac{j^{*,\pm}(\d s)}{\sigma(\phi^{-1}(s,v))} 
\quad \text{as $ n \to \infty $}. 
\label{eq: regvar jv2}
\end{align}
Suppose, moreover, that $ (0,j^*,0) $ satisfies {\bf (N0)}-{\bf (N3)}. 
Then the same assertions as Theorem \ref{hom2} hold. 
\end{Cor}

\Proof[Proof of Corollary \ref{SG2}]{
Let $ S'' = \{ +,- \} \times [0,1] \times (0,\infty ) $ 
and let $ \mu(\{ \pm \} \times \d v \d r) = \d r $. 
For $ v \in [0,1] $, we define 
a function $ J^{0,\pm}_v:(0,\infty ) \to (0,\infty ) $ by 
\begin{align}
J^{0,\pm}_vy = \inf \cbra{ r>0 : \int_{(0,r]} \frac{j^{\pm}_v(\d s)}{\sigma(\phi^{-1}(s,v))} > y } 
\label{}
\end{align}
and define a function $ J:\{ +,- \} \times [0,1] \times (0,\infty ) \to G \setminus \{ 0 \} $ 
by $ J(\pm,v,y) = (J^{0,\pm}_vy)v $. 
In the same way as $ J $ we define $ J^* $ 
with $ j^{\pm}_v $ being replaced by $ j^{*,\pm}_v $. 
For $ n \in \bN $, set 
$ J_n(\pm,v,y) = c^{-\alpha n} J(\pm,v,c^{\alpha (\kappa -\beta) n}y) $. 
By the assumption \eqref{eq: regvar jv2}, we see that 
{\bf (C5)} and {\bf (C6)} are satisfied. 
We can thus apply Theorem \ref{hom2}. 
}

We may expect the following. 

\begin{Conj} \label{conj}
For $ \{ x_n \} \subset G $ such that $ x_n \to 0 $, 
it holds that 
$ T_{x_n} \to 0 $, $ \vn $-a.e. 
\end{Conj}

We do not know whether Conjecture \ref{conj} is true or not. 
If Conjecture \ref{conj} is true, 
then we can easily obtain the following. 

\begin{Conj}[jumping-in vanishing case]
Let $ (\rho,j,\varsigma) $ satisfy {\bf (CG1)}-{\bf (CG4)}. 
Suppose, moreover, that 
$ \int_{G \setminus \{ 0 \}} \frac{j(\d x)}{\sigma(x)} < \infty $. Set 
\begin{align}
\rho^*_{\pm} 
:= \rho_{\pm} 
+ \int_{G_{\pm} \setminus \{ 0 \}} \frac{j(\d x)}{\sigma(x)} . 
\label{}
\end{align}
Then the same assertions as Theorem \ref{hom} hold. 
\end{Conj}

\section{Proof of the general limit theorem} \label{sec: pr glt}

\subsection{Piecing proposition} \label{sec: piece}

Let us prove Proposition \ref{piece}. 

\Proof[Proof of Proposition \ref{piece}]{
Let us write $ \eta(l) $, $ L(t) $ and $ X(t) $ 
simply for $ \eta(p,\varsigma;t) $, $ L(p,\varsigma;t) $ and $ X(p,\varsigma;t) $. 

We prove \eqref{eq: occup}. 
For $ t \ge 0 $, we have 
\begin{align}
\int_0^t 1_{\{ X(s)=0 \}} \d s 
=& t - \int_0^t 1_{\{ X(s) \neq 0 \}} \d s 
\label{} \\
=& t - \sum_{l<L(t)} \int_{\eta(l-)}^{\eta(l)} 1_{\{ X(s) \neq 0 \}} \d s 
- \int_{\eta(L(t)-)}^t 1_{\{ X(s) \neq 0 \}} \d s 
\label{} \\
=& t - \sum_{l<L(t)} T_0(p(l)) 
- \{ t - \eta(L(t)-) \} , 
\label{}
\end{align}
which is equal to $ \varsigma L(t) $ by the definition \eqref{eq: eta}. 
Thus we obtain \eqref{eq: occup}. 

Let us assume that the condition {\bf (N3)} is satisfied. 
By {\bf (N3)}, we see that, for any $ n \in \bN $, 
there are at most finitely many $ l \le n $ such that $ \| p^{(l)} \| \ge 1/n $. 
This shows that, 
if there exists a sequence $ l_n $ converging to $ l $ such that $ p^{(l_n)} \neq o $, 
then it implies that $ \| p^{(l_n)} \| \to 0 $. 
We now let $ t_0 \ge 0 $ and we prove that $ X(t) $ is c\`adl\`ag at $ t=t_0 $. 
Set $ l_0 = L(t_0) $, $ t_1 = \eta(l_0-) $ and $ t_2 = \eta(l_0) $ 
so that $ t_1 \le t_0 \le t_2 $. 
We divide the proof into three cases. 
\\ (i) 
Suppose that $ t_1 \le t_0 < t_2 $. 
We then have $ X(t) = p^{(l_0)}(t-t_1) $ for all $ t_1 \le t < t_2 $. 
This shows that $ X(t) $ is right continuous at $ t=t_0 $ 
and has left limit at $ t=t_0 $ except when $ t_1 = t_0 $. 
If $ t_1 = t_0 $, 
we see, by the above remark, that $ X(t_0-)=0 $. 
\\ (ii) 
Suppose that $ t_1 < t_0 = t_2 $. 
We then have $ X(t) = p^{(l_0)}(t-t_1) $ for all $ t_1 \le t < t_0 $. 
This shows that $ X(t) $ has left limit at $ t=t_0 $. 
Since $ \eta $ is strictly increasing, 
there is no $ l $ such that $ \eta(l_0)=\eta(l-) $, 
and hence $ X(t_0)=0 $ by definition of $ X $. 
If there exists a sequence $ l_n $ decreasing to $ l_0 $ 
such that $ p^{(l_n)} \neq o $, we have $ \| p^{(l_n)} \| \to 0 $ by the above remark, 
and hence we obtain $ X(t_0+)=0 $. 
Otherwise, we have $ X(t_0+)=0 $ by definition of $ X $. 
\\ (iii) 
Suppose that $ t_1 = t_0 = t_2 $. 
We then easily see that $ X(t_0-) = X(t_0) = X(t_0+) = 0 $.

Let us assume that $ \vn(D)=\infty $ but that the condition {\bf (N3)} is not satisfied. 
We then have $ \vn \{ w \in D : \| w \| \ge r_0 \} = \infty $ for some $ r_0>0 $. 
Hence there exists a sequence $ l_n $ decreasing to 0 such that 
$ \| p^{(l_n)} \| \ge r_0 $ for all $ n $. 
By {\bf (N1)}, we have $ T_0(p^{(l_n)}) \to 0 $, 
and hence we obtain $ \limsup_{t \to 0+} \| X(t) \| \ge r_0 $. 
Since we have $ X(0) = 0 $ by the definition of $ X $, 
we see that $ X $ is not right-continuous. 

Therefore we conclude the proof. 
}

\subsection{Useful lemmas}

For the proof of Theorem \ref{glt}, we need two lemmas. 
The first one is the following. 

\begin{Lem} \label{lem1}
Let $ w_1,w_2,\ldots,w_{\infty } \in D^0 $. 
Suppose that $ w_n \to w_{\infty } $ in $ D $ 
and that $ T_0(w_n) \to T_0(w_{\infty }) $. 
Then one has $ \| w_n \| \to \| w_{\infty } \| $. 
\end{Lem}

\Proof{
By the assumption that $ w_n \to w_{\infty } $ in $ D $, 
we may take transformations $ I_1,I_2,\ldots,I_{\infty } $ of $ [0,\infty ) $ 
such that $ I_n \to I $ uc and $ w_n - w_{\infty } \circ I_n \to 0 $ uc. 
Let $ \eps>0 $. Then we may choose $ N \in \bN $ so that 
for any $ n \ge N $ we have 
\begin{align}
\sup_{t \le T_0(w_{\infty })+1} |w_n(t) - w_{\infty }(I_n(t)) | < \eps 
\label{}
\end{align}
and we have 
\begin{align}
T_0(w_n) \le T_0(w_{\infty }) + 1 
\quad \text{and} \quad 
\sup_{t \le T_0(w_{\infty }) + 1} |I_n(t)-t| \le 1 . 
\label{}
\end{align}
For $ n \ge N $, we have 
\begin{align}
\| w_n \| 
=& \sup_{t \le T_0(w_n)} |w_n(t)| 
\label{} \\
\le& \sup_{t \le T_0(w_{\infty })+1} |w_{\infty }(I_n(t))| 
+ \sup_{t \le T_0(w_{\infty })+1} | w_n(t) - w_{\infty }(I_n(t)) | 
\label{} \\
\le& 
\| w_{\infty } \| + \eps , 
\label{}
\end{align}
and also we have 
\begin{align}
\| w_{\infty } \| 
=& \sup_{t \le T_0(w_{\infty })} |w_{\infty }(t)| 
\label{} \\
=& \sup_{t \le T_0(w_{\infty })+1} |w_{\infty }(I_n(t))| 
\label{} \\
\le& \sup_{t \le T_0(w_{\infty })+1} |w_n(t)| 
+ \sup_{t \le T_0(w_{\infty })+1} | w_n(t) - w_{\infty }(I_n(t)) | 
\label{} \\
\le& 
\| w_n \| + \eps . 
\label{}
\end{align}
Hence we obtain $ |\| w_n \| - \| w_{\infty } \|| < \eps $. 
The proof is complete. 
}

\begin{Rem}
We cannot remove the assumption $ T_0(w_n) \to T_0(w_{\infty }) $ 
from Lemma \ref{lem1}. 
In fact, if we set 
\begin{align}
w_n(t) = 
\begin{cases}
1 & \text{if $ 0 \le t < 1 $}, \\
1/n & \text{if $ 1 \le t < n $}, \\
2 & \text{if $ n \le t < n+1 $}, \\
0 & \text{if $ t \ge n+1 $}, 
\end{cases}
\quad 
w_{\infty }(t) = 
\begin{cases}
1 & \text{if $ 0 \le t < 1 $}, \\
0 & \text{if $ t \ge 1 $}, 
\end{cases}
\label{}
\end{align}
then we have $ w_n \to w_{\infty } $ in $ D $ 
but $ \| w_n \| = 2 $ and $ \| w_{\infty } \| = 1 $. 
\end{Rem}

The following lemma is partly taken from Bartle \cite{MR600921}. 

\begin{Lem}[Bartle \cite{MR600921}] \label{lem: a.u.}
Let $ f_1,f_2,\ldots,f_{\infty } $ be real-valued functions 
defined on a measure space $ (E,\cE,\nu) $. 
Then the following statements are equivalent: 
\begin{enumerate}
\item 
$ f_n \to f_{\infty } $, $ \nu $-almost uniformly; 
\item 
$ \nu \rbra{ \bigcup_{n=N}^{\infty } \{ |f_n-f_{\infty }| \ge r \} \Big. } \to 0 $ 
as $ N \to \infty $ for all $ r>0 $; 
\item 
$ f_n \to f_{\infty } $, $ \nu $-a.e. 
and for any $ r>0 $ 
there exists $ N $ such that 
\\
$ \nu \rbra{ \bigcup_{n=N}^{\infty } \{ |f_n-f_{\infty }| \ge r \} \Big. } < \infty $. 
\end{enumerate}
\end{Lem}

The proof can be found in Bartle \cite[Theorem 1.7]{MR600921}, 
and so we omit it.

\begin{Lem} \label{lem2}
Let $ \nu $ be a $ \sigma $-finite measure on a Polish space $ E $ 
and let $ \Phi_1,\Phi_2,\ldots,\Phi_{\infty } $ be measurable mappings 
from $ E $ to $ D $ such that 
we have $ \Phi_n \to \Phi_{\infty } $ in $ D $, $ \nu $-a.e. 
and $ T_0(\Phi_n) \to T_0(\Phi_{\infty }) $, $ \nu $-a.e. 
Suppose that $ \nu \rbra{ \| \Phi_n \| \ge r } < \infty $ 
for all $ n=1,2,\ldots,\infty $ and all $ r>0 $. 
Then the following statements are equivalent: 
\begin{enumerate}
\item 
$ \| \Phi_n \| \to \| \Phi_{\infty } \| $, $ \nu $-almost uniformly; 
\item 
$ \nu \rbra{ \bigcup_{n=1}^{\infty } \cbra{ \| \Phi_n \| \ge r } \Big. } < \infty $ 
for all $ r>0 $. 
\end{enumerate}
\end{Lem}

\Proof{
For $ r>0 $ and $ n \in \bN \cup \{ \infty \} $, we write 
\begin{align}
A^r_n = \cbra{ \| \Phi_n \| \ge r } 
\quad \text{and} \quad 
B^r_n = \cbra{ | \| \Phi_n \| - \| \Phi_{\infty } \| | \ge r } . 
\label{}
\end{align}

Suppose (i) is satisfied. 
Let $ r>0 $ be fixed. 
We then see that 
$ \nu \rbra{ \bigcup_{n=N}^{\infty } B^{r/2}_n } < \infty $ for some $ N $ 
by Lemma \ref{lem: a.u.}. 
Since we have $ \nu \rbra{ \bigcup_{n=1}^N B^{r/2}_n } < \infty $ by the assumption, 
we obtain $ \nu \rbra{ \bigcup_{n=1}^{\infty } B^{r/2}_n } < \infty $. 
Since we have 
$ \bigcup_{n=1}^{\infty } A^r_n \subset 
\rbra{ \bigcup_{n=1}^{\infty } B^{r/2}_n } \cup A^{r/2}_{\infty } $, 
we see that (ii) is satisfied. 

Suppose (ii) is satisfied. 
Let $ r>0 $ be fixed. 
We then see that 
$ \nu \rbra{ \bigcup_{n=1}^{\infty } B^r_n } 
\le \nu \rbra{ \bigcup_{n=1}^{\infty } A^{r/2}_n } 
+ \nu \rbra{ A^{r/2}_{\infty } \Big. } < \infty $. 
Note that, by Lemma \ref{lem1}, 
we have $ \| \Phi_n \| \to \| \Phi_{\infty } \| $, $ \nu $-a.e. 
Hence, by Lemma \ref{lem: a.u.}, 
we see that (i) is satisfied. 
}

\subsection{General limit theorem} \label{sec: prf glt}

We now proceed to the proof of Theorem \ref{glt}. 

\Proof[Proof of Theorem \ref{glt}]{
Let $ E $, $ \nu $, $ \Phi_1,\Phi_2,\ldots $ and $ \Phi_{\infty } $ 
be as in Definition \ref{def: conv}. 
Taking a subsequence if necessary, 
we may take $ N=1 $ in Condition {\bf (G5)} of Definition \ref{def: conv}. 

Let $ p = (p^{(l)})_{l \in D(p)} $ 
be a Poisson point process on $ E $ with characteristic measure $ \nu $. 
For $ n \in \bN \cup \{ \infty \} $, we define 
\begin{align}
q_n^{(l)} = 
\begin{cases}
\Phi_n(p^{(l)}) 
& \text{if $ l \in D(p) $}, \\
o & \text{otherwise}. 
\end{cases}
\label{}
\end{align}
We then see that 
$ q_n = (q_n^{(l)})_{l \ge 0} $ is a Poisson point process on $ D $ outside $ o $ 
with characteristic measure $ (\nu \circ \Phi_n^{-1}) |_{D \setminus \{ o \}} $, 
which is equal to $ \vn_n $ 
by Assumption {\bf (A2)} and Condition {\bf (G1)}. 
Hence $ q_n $ is a realization of $ p_n $, 
so that we may assume without loss of generality that $ p_n = q_n $. 

For $ m \in \bN $, we set 
\begin{align}
\Lambda_m = \{ l \in [0,m] : \| p_n^{(l)} \| \ge 1/m 
\ \text{for some $ n \in \bN \cup \{ \infty \} $} \} . 
\label{}
\end{align}
By Assumption {\bf (A2)} and Condition {\bf (G3)$ ' $}, we see that 
$ \sharp \Lambda_m < \infty $ a.s. for all $ m \in \bN $. 
By this fact 
and by the assumptions {\bf (A1)}-{\bf (A3)}, 
we see that there exists an event $ \Omega^* $ of probability one 
such that for any sample point belonging to $ \Omega^* $ we have the following: 
\begin{description}
\item[(L1)] 
for any $ l \ge 0 $, $ p_n^{(l)} \to p_{\infty }^{(l)} $ in $ D $; 
\item[(L2)] 
for any $ l \ge 0 $, $ T_0(p_n^{(l)}) \to T_0(p_{\infty }^{(l)}) $ in $ [0,\infty ] $; 
\item[(L3)] 
for any $ l \ge 0 $, $ \displaystyle \tau(l) := \sup_{n \ge 1} T_0(p_n^{(l)}) $ satisfies 
$ \displaystyle \sum_{s \le l} \tau(s) < \infty $; 
\item[(L4)] 
$ p^{(l)}_{\infty }(T_0(p_{\infty }^{(l)})-) = 0 $; 
\item[(L5)] 
for any $ m \in \bN $, $ \sharp \Lambda_m < \infty $. 
\end{description}
In what follows we pick and fix a sample point belonging to $ \Omega^* $. 

Since we have 
\begin{align}
\eta_n(l) = \varsigma_n l + \sum_{s \le l} T_0(p_n^{(l)}) , 
\label{}
\end{align}
we see that, for any $ l_0>0 $, 
\begin{align}
G_n(l_0) :=& 
\sup_{l \le l_0} |\eta_n(l)-\eta_{\infty }(l)| 
\label{} \\
\le& |\varsigma_n-\varsigma_{\infty }| l_0 + \sum_{s \le l_0} |T_0(p_n^{(l)}) - T_0(p_{\infty }^{(l)})| . 
\label{}
\end{align}
By {\bf (L2)}-{\bf (L3)} and by {\bf (A4)}, we apply the dominated convergence theorem 
to see that 
\begin{align}
G_n(l_0) \to 0 \quad \text{as $ n \to \infty $ for all $ l_0>0 $}. 
\label{eq: Gn}
\end{align}
Hence we obtain $ \eta_n \to \eta_{\infty } $ in $ D $. 
Moreover, since $ \eta_{\infty } $ is strictly increasing, 
we use \cite[Theorem 7.2]{MR561155} to obtain 
\begin{align}
H_n(t_0) := \sup_{t \le t_0} |L_n(t)-L_{\infty }(t)| \to 0 
\label{eq: Hn}
\end{align}
as $ n \to \infty $ for all $ t_0>0 $. 
Hence we obtain $ L_n \to L_{\infty } $ in $ D $.

It remains to prove that $ X_n \to X_{\infty } $ in $ D $. 
So we take an arbitrary subsequence 
and denote it by the same symbol as the original sequence. 
It suffices to prove that 
we can extract a further subsequence along which $ X_n \to X_{\infty } $ in $ D $. 
We divide the proof into several steps. 

{\bf Step 1.} 
For $ l \ge 0 $, 
since $ p_n^{(l)} \to p_{\infty }^{(l)} $ in $ D $, we see that 
there exist transformations $ I^{(l)}_1,I^{(l)}_2,\ldots,I^{(l)}_{\infty } $ of $ [0,\infty ) $ 
such that each $ I^{(l)}_n $ is bijective, continuous and increasing 
and that we have 
\begin{align}
F^{(l)}_n(t_0) := 
\sup_{t \le t_0} \absol{ I^{(l)}_n(t)-t } 
+ 
\sup_{t \le t_0} 
\absol{ p^{(l)}_n(t) - p^{(l)}_{\infty }(I^{(l)}_n(t)) \Big. } 
\tend{}{n \to \infty } 0 
\label{eq: Fn}
\end{align}
for all $ t_0>0 $. 
For $ n \in \bN $ and $ m \in \bN $, we set 
\begin{align}
F_n(l) = F^{(l)}_n(T_0(p_{\infty }^{(l)}) + 3) 
\quad \text{for $ l \ge 0 $} 
\label{}
\end{align}
and set 
\begin{align}
M_{n,m} = \max_{l \in \Lambda_m} 
\cbra{ F_n(l) + G_n(l) + \absol{ T_0(p_n^{(l)}) - T_0(p_{\infty }^{(l)}) \Big. } } . 
\label{}
\end{align}
By \eqref{eq: Fn}, \eqref{eq: Gn}, {\bf (L2)} and {\bf (L5)}, we see that 
$ M_{n,m} \to 0 $ as $ n \to \infty $ for all fixed $ m \in \bN $. 
Thus we may take a subsequence $ \{ n_1(n) \}_{n \in \bN} $ such that 
$ M_{n_1(n),n} < 1/n $ for all $ n \in \bN $. 
Writing $ p_n $ simply for $ p_{n_1(n)} $, 
we may assume without loss of generality that 
$ M_n := M_{n,n} < 1/n $ for all $ n \in \bN $. 

{\bf Step 2.} 
We construct a transformation $ I_n $ of $ [0,\infty ) $. 

We modify the transformation $ I^{(l)}_n $ around $ t=T_0(p_n^{(l)}) $. 
Note that 
\begin{align}
|T_0(p_n^{(l)}) - T_0(p_{\infty }^{(l)}) | \le M_n < 1/n 
\label{}
\end{align}
and that 
\begin{align}
T_0(p_n^{(l)}) + 2/n < T_0(p_{\infty }^{(l)}) + 3/n \le T_0(p_{\infty }^{(l)}) + 3 . 
\label{}
\end{align}
Since $ F_n(l) \le M_n < 1/n $, we have 
\begin{align}
s^{(l)}_n :=& I^{(l)}_n(T_0(p_n^{(l)})-2/n) \le (T_0(p_n^{(l)})-2/n)+M_n < T_0(p_{\infty }^{(l)}) , 
\label{} \\
t^{(l)}_n :=& I^{(l)}_n(T_0(p_n^{(l)})+2/n) \ge (T_0(p_n^{(l)})+2/n)-M_n > T_0(p_{\infty }^{(l)}) . 
\label{}
\end{align}
We define 
\begin{align}
\tilde{I}^{(l)}_n(t) = 
\begin{cases}
0 & \text{if $ t=0 $}, \\
I^{(l)}_n(t) & \text{if $ |t-T_0(p_n^{(l)})| \ge 2/n $}, \\
T_0(p_{\infty }^{(l)}) & \text{if $ t=T_0(p_n^{(l)}) $}, \\
{\rm linear} & \text{otherwise}. 
\end{cases}
\label{}
\end{align}
Since $ s^{(l)}_n < T_0(p_{\infty }^{(l)}) < t^{(l)}_n $, 
we see that $ \tilde{I}^{(l)}_n $ is well-defined, bijective, increasing and continuous 
and that 
\begin{align}
\tilde{I}^{(l)}_n \rbra{ T_0(p_n^{(l)}) \Big. } = T_0(p_{\infty }^{(l)}) . 
\label{}
\end{align}
We note that 
\begin{align}
\sup_{t \ge 0} \absol{ \tilde{I}^{(l)}_n(t) - I^{(l)}_n(t) } 
=& \absol{ T_0(p_{\infty }^{(l)}) - I^{(l)}_n(T_0(p_n^{(l)})) \Big. } 
\label{} \\
\le& \absol{ T_0(p_{\infty }^{(l)}) - T_0(p_n^{(l)}) \Big. } + M_n 
< 2/n . 
\label{}
\end{align}

Denote $ l_n = \max \Lambda_n $ if $ \Lambda_n \neq \emptyset $ 
and $ l_n = 0 $ if $ \Lambda_n = \emptyset $. 
We now define 
\begin{align}
I_n(t) = 
\begin{cases}
0 & \text{if $ t=0 $}, \\
\eta_{\infty }(l-) + \tilde{I}^{(l)}_n(t-\eta_n(l-)) 
& \text{if $ \eta_n(l-) \le t < \eta_n(l) $ for some $ l \in \Lambda_n $} , \\
t - \eta_n(l_n) + \eta_{\infty }(l_n) 
& \text{if $ t \ge \eta_n(l_n) $}, \\
{\rm linear} & \text{otherwise}. 
\end{cases}
\label{eq: In}
\end{align}
We then see that $ I_n $ is bijective, increasing and continuous. 

{\bf Step 3.} 
We prove that $ I_n \to I $ uc. 
Let $ t_0>0 $ be fixed. 
Since $ L_n(t_0) \to L_{\infty }(t_0) $, 
there exists $ N $ such that 
$ L_n(t_0) \le L_{\infty }(t_0) + 1 $ for all $ n \ge N $. 

If $ t $ is such that $ \eta_n(l-) \le t \le \eta_n(l) $ for some $ l \in \Lambda_n $, 
we have, for $ n \ge N $, 
\begin{align}
|I_n(t)-t| 
\le& |\eta_{\infty }(l-)-\eta_n(l-)| 
+ \absol{ \tilde{I}^{(l)}_n(t-\eta_n(l-)) - (t-\eta_n(l-)) } 
\label{} \\
\le& \sup_{u \le L_n(t_0)} |\eta_{\infty }(u)-\eta_n(u)| 
+ \sup_{s \le T_0(p_n^{(l)})} \absol{ I^{(l)}_n(s)-s } + 2/n 
\label{} \\
\le& G_n(L_{\infty }(t_0)+1) 
+ \sup_{s \le T_0(p_{\infty }^{(l)})+1/n} \absol{ I^{(l)}_n(s)-s } + 2/n 
\label{} \\
\le& G_n(L_{\infty }(t_0)+1) + 3/n . 
\label{}
\end{align}
Otherwise, we have, by linearity, 
\begin{align}
|I_n(t)-t| 
\le 
\sup_{u \le L_n(t_0)} |\eta_{\infty }(u)-\eta_n(u)| 
\le G_n(L_{\infty }(t_0)+1) . 
\label{}
\end{align}
Thus we obtain 
\begin{align}
\sup_{t \le t_0} |I_n(t)-t| 
\le G_n(L_{\infty }(t_0)+1) + 3/n 
\tend{}{n \to \infty } 0 . 
\label{}
\end{align}
This shows that $ I_n \to I $ uc.

{\bf Step 4.} 
Let $ m \in \bN $ be fixed. 
For $ n \ge m $ and for $ l \in \Lambda_m $ ($ \subset \Lambda_n $), 
we estimate the supremum over $ t \ge 0 $ of 
\begin{align}
K^{(l)}_n(t) := \absol{ p^{(l)}_n(t) - p^{(l)}_{\infty }(\tilde{I}^{(l)}_n(t)) } . 
\label{}
\end{align}

If $ t \le T_0(p_n^{(l)})-2/n $, 
we have $ t \le T_0(p_{\infty }^{(l)})+3 $, and hence we have 
\begin{align}
K^{(l)}_n(t) 
= \absol{ p^{(l)}_n(t) - p^{(l)}_{\infty }(I^{(l)}_n(t)) } 
\le M_n < 1/n . 
\label{}
\end{align}
If $ T_0(p_n^{(l)})-2/n < t < T_0(p_n^{(l)}) $, 
we have $ t < T_0(p_{\infty }^{(l)})+1/n $ and 
\begin{align}
\tilde{I}^{(l)}_n(t) 
\ge I^{(l)}_n(t) - 2/n 
\ge t - 3/n 
\ge T_0(p_{\infty }^{(l)}) - 4/n , 
\label{}
\end{align}
which yields 
\begin{align}
K^{(l)}_n(t) 
\le& 
\absol{ p^{(l)}_n(t) - p^{(l)}_{\infty }(I^{(l)}_n(t)) } 
+ | p^{(l)}_{\infty }(I^{(l)}_n(t)) | 
+ | p^{(l)}_{\infty }(\tilde{I}^{(l)}_n(t)) | 
\label{} \\
\le& 
1/n + 2 \sup_{s \ge T_0(p_{\infty }^{(l)})-4/n} | p^{(l)}_{\infty }(s) | . 
\label{}
\end{align}
If $ t \ge T_0(p_n^{(l)}) $, 
we have $ p^{(l)}_n(t) = p^{(l)}_{\infty }(I^{(l)}_n(t)) = 0 $, 
so that $ K^{(l)}_n(t) = 0 $. 

Therefore we obtain 
\begin{align}
\sup_{t \ge 0} K^{(l)}_n(t) 
\le 
1/n + 2 \sup_{s \ge T_0(p_{\infty }^{(l)})-4/n} | p^{(l)}_{\infty }(s) | , 
\label{}
\end{align}
which converges to 0 by {\bf (L4)}. 

We now set 
\begin{align}
\tilde{F}_n(l) = \sup_{t \le T_0(p_{\infty }^{(l)})+3} K^{(l)}_n(t) 
\label{}
\end{align}
and 
\begin{align}
\tilde{M}_{n,m} = \max_{l \in \Lambda_m} \tilde{F}_n(l) . 
\label{}
\end{align}
We then have $ \tilde{M}_{n,m} \to 0 $ as $ n \to \infty $ 
for all fixed $ m \in \bN $. 
Hence we may take a subsequence $ \{ n_2(n) \}_{n \in \bN} $ 
such that $ \tilde{M}_{n_2(n),n}<1/n $ for all $ n \in \bN $. 
Writing $ p_n $ simply for $ p_{n_2(n)} $, 
we may assume without loss of generality that 
$ \tilde{M}_n := \tilde{M}_{n,n} < 1/n $ for all $ n \in \bN $.

{\bf Step 5.} 
Let $ t_0>0 $. We estimate the supremum over $ 0 \le t \le t_0 $ of 
\begin{align}
K_n(t) := 
\absol{ X_n(t) - X_{\infty }(I_n(t)) } . 
\label{}
\end{align}

Taking a subsequence if necessary, 
we may assume without loss of generality that 
\begin{align}
\sup_{t \le t_0} |I_n(t)-t| < 1/n 
\quad \text{for all $ n \in \bN $}. 
\label{}
\end{align}
Let $ N $ be such that $ L_n(t_0) \le L_{\infty }(t_0)+1 $ for $ n \ge N $. 
We take $ N $ large enough to satisfy $ N > L_{\infty }(t_0+1)+1 $. 
Let $ n \ge N $. 

\noindent
\un{Case 1}: $ \eta_n(l-) \le t < \eta_n(l) $ for some $ l \in \Lambda_n $. 
Since $ \eta_{\infty }(l-) \le I_n(t) < \eta_{\infty }(l) $, we have 
\begin{align}
K_n(t) 
= \absol{ p^{(l)}_n(t-\eta_n(l-)) - p^{(l)}_{\infty }(\tilde{I}^{(l)}_n(t-\eta_n(l-)) } 
\le \tilde{M}_n < 1/n . 
\label{}
\end{align}
\un{Case 2}: $ \eta_n(l-) \le t < \eta_n(l) $ for some $ l \notin \Lambda_n $. 
Since $ n > L_{\infty }(t_0)+1 \ge L_n(t_0) $, we have 
$ \eta_n(n) > t_0 \ge t \ge \eta_n(l-) $. 
Hence we have $ n>l $. 
Since $ l \notin \Lambda_n $, we have $ \| p_n^{(l)} \| < 1/n $. 
We now obtain 
\begin{align}
K_n(t) 
=& \absol{ p^{(l)}_n(t-\eta_n(l-)) - p^{(l)}_{\infty }(\tilde{I}^{(l)}_n(t-\eta_n(l-)) } 
\label{} \\
\le& \| p_n^{(l)} \| + \| p_{\infty }^{(l)} \| 
\le \tilde{M}_n < 2/n . 
\label{}
\end{align}
\un{Case 3}: there is no $ l \ge 0 $ such that $ \eta_n(l-) \le t < \eta_n(l) $. 
In this case we have $ X_n(t)=0 $. 
We divide the proof into three subcases. 
\\
\un{Case 3-1}: 
there is no $ l \ge 0 $ such that $ \eta_{\infty }(l-) \le I_n(t) < \eta_{\infty }(l) $. 
In this case we have $ X_{\infty }(I_n(t))=0 $, 
so that $ K_n(t)=0 $. 
\\
\un{Case 3-2}: 
$ \eta_{\infty }(l-) \le I_n(t) < \eta_{\infty }(l) $ for some $ l \ge 0 $ 
with $ \| p_{\infty }^{(l)} \| < 1/n $. 
In this case we have $ | X_{\infty }(I_n(t)) | < 1/n $, 
so that $ K_n(t) < 1/n $. 
\\
\un{Case 3-3}: 
$ \eta_{\infty }(l-) \le I_n(t) < \eta_{\infty }(l) $ for some $ l \ge 0 $ 
with $ \| p_{\infty }^{(l)} \| \ge 1/n $. 
In this case, we have 
\begin{align}
l \le L_{\infty }(I_n(t)) \le L_{\infty }(t_0+1) < n , 
\label{}
\end{align}
so that we have $ l \in \Lambda_n $. 
Thus we obtain 
\begin{align}
\sup_{t \le t_0} \absol{ p^{(l)}_n(t) - p^{(l)}_{\infty }(I^{(l)}_n(t)) \Big. } 
\le M_n < 1/n , 
\label{}
\end{align}
so that 
\begin{align}
\| p_n^{(l)} \| \ge \| p_{\infty }^{(l)} \| 
- \sup_{t \le t_0} \absol{ p^{(l)}_n(t) - p^{(l)}_{\infty }(I^{(l)}_n(t)) \Big. } 
> 1/n - 1/n = 0 . 
\label{}
\end{align}
From this we see that $ \eta_n(l-) < \eta_n(l) $ and 
by the definition \eqref{eq: In} of $ I_n $ we see that 
$ \eta_n(l-) \le s < \eta_n(l) $ implies $ \eta_{\infty }(l-) \le I_n(s) < \eta_{\infty }(l) $. 
Thus we obtain $ \eta_n(l-) \le t < \eta_n(l) $, 
which is a contradiction. 

From all the arguments above, we obtain 
\begin{align}
\sup_{t \le t_0} K_n(t) \le 2/n . 
\label{}
\end{align}
This shows that $ X_n - X_{\infty } \circ I_n \to 0 $ uc. 

The proof is therefore complete. 
}

\section{Proof of the homogenization theorem} \label{sec: pr hom} 

\subsection{Scaling property for local time and excursion measure} \label{sec: prsca}

Let us prove Lemma \ref{stag0}. 

\Proof[Proof of Lemma \ref{stag0}]{
We deal with processes under $ \bP_0 $. 
Using {\bf (S2)} and \eqref{eq: stag}, we have 
\begin{align}
c^{-\alpha \kappa n} \int_0^{c^n t} 1_{\{ X(s)=0 \}} \d s 
\law \int_0^t 1_{\{ X(s)=0 \}} \d s . 
\label{eq: by S2}
\end{align}
Hence, using {\bf (S1)} and then using \eqref{eq: by S2}, we obtain 
\begin{align}
\int_0^t 1_{\{ X(s)=0 \}} \d s 
\law& \int_0^t 1_{\{ (\Psi_{\alpha }^n X)(s)=0 \}} \d s 
\label{} \\
=& c^{-n} \int_0^{c^n t} 1_{\{ X(s)=0 \}} \d s 
\label{} \\
\law& c^{-(1-\alpha \kappa )n} \int_0^t 1_{\{ X(s)=0 \}} \d s . 
\label{}
\end{align}
Since $ 0<\alpha \kappa <1 $ by {\bf (S0)}, 
the last quantity converges in law to 0 as $ n \to \infty $. 
}

To prove Proposition \ref{excscale} and for the later use, 
we prove the following lemma. 

\begin{Lem} \label{scaling}
Let $ (\vn,\varsigma) $ be the pair consisting of a $ \sigma $-finite measure on $ D $ 
and a non-negative constant $ \varsigma $ 
and suppose that $ (\vn,\varsigma) $ satisfies Conditions {\bf (N0)}-{\bf (N3)}. 
Let $ c $, $ \alpha $ and $ \gamma $ be positive constants. 
Let $ \{ p,\bP \} $ be a Poisson point process on $ D $ outside $ o $ 
with characteristic measure $ \vn $ and denote 
\begin{align}
\tilde{p}^{(l)} = \Psi_{\alpha } p^{(c^{\gamma } l)} 
, \quad 
\tilde{\vn} = c^{\gamma } \vn \circ \Psi_{\alpha }^{-1} 
, \quad 
\tilde{\varsigma} = c^{-(1-\gamma )} \varsigma . 
\label{}
\end{align}
Then it holds 
that the pair $ (\tilde{\vn},\tilde{\varsigma}) $ 
satisfies {\bf (N0)}-{\bf (N3)}, 
that $ \{ \tilde{p},\bP \} $ is a Poisson point process on $ D $ outside $ o $ 
with characteristic measure $ \tilde{\vn} $ 
and that 
\begin{align}
X(\tilde{p},\tilde{\varsigma}) = \Psi_{\alpha } X(p,\varsigma) 
, \quad 
L(\tilde{p},\tilde{\varsigma}) = \Psi_{\gamma } L(p,\varsigma) 
, \quad 
\eta(\tilde{p},\tilde{\varsigma}) = \hat{\Psi}_{\gamma } \eta(p,\varsigma) . 
\label{eq: tileta}
\end{align}
\end{Lem}

\Proof{
For $ l \ge 0 $ and $ A \in \cB(D \setminus \{ o \}) $, we have 
\begin{align}
\sharp \{ s \le l : \tilde{p}^{(s)} \in A \} 
= \sharp \{ s \le l : p^{(c^{\gamma} s)} \in \Psi_{\alpha }^{-1} A \} 
= N_{c^{\gamma} l}(\Psi_{\alpha }^{-1} A) , 
\label{}
\end{align}
where we note that $ \Psi_{\alpha }^{-1} \{ o \} = \{ o \} $. 
It is now obvious that $ \tilde{p} $ is again 
a Poisson point process on $ D $ outside $ o $ 
whose characteristic measure is 
equal to $ \tilde{\vn} = c^{\gamma} \vn \circ \Psi_{\alpha }^{-1} $. 
Hence we obtain 
\begin{align}
\hat{\Psi}_{\gamma} \eta(p,\varsigma;l) 
=& c^{-1} \cbra{ \varsigma c^{\gamma} l + \sum_{s \le c^{\gamma} l} T_0(p^{(s)}) } 
\label{} \\
=& \tilde{\varsigma} l + \sum_{s \le l} T_0(p^{(c^{\gamma} s)}) / c 
\label{} \\
=& \tilde{\varsigma} l + \sum_{s \le l} T_0(\tilde{p}^{(s)}) , 
\label{}
\end{align}
which yields 
$ \eta(\tilde{p},\tilde{\varsigma}) = \hat{\Psi}_{\gamma} \eta(p,\varsigma) $. 
The other identities of \eqref{eq: tileta} are now obvious. 
}

We now prove Proposition \ref{excscale}. 

\Proof[Proof of Proposition \ref{excscale}]{
Suppose that {\bf (S2)$ ' $} is satisfied. 
Let $ p $ be a Poisson point process on $ D $ outside $ o $ 
with characteristic measure $ \vn $. 
Denote $ \tilde{p}^{(l)} = \Psi_{\alpha } p^{(c^{\alpha \kappa } l)} $. 
By Lemma \ref{scaling} and {\bf (S2)$ ' $}, 
we see that $ \tilde{p} $ is equal in law to $ p $. 
On one hand, we see that 
\begin{align}
\{ X(\tilde{p},0), L(\tilde{p},0) \} 
\law 
\{ X(p,0) , L(p,0) \} . 
\label{}
\end{align}
On the other hand, we see that 
\begin{align}
L(\tilde{p},0;t) 
=& \inf \cbra{ l \ge 0 : \sum_{s \le l} T_0(\Psi_{\alpha } p^{(c^{\alpha \kappa } s)}) > t } 
\label{} \\
=& \inf \cbra{ l \ge 0 : \sum_{s \le l} T_0(p^{(c^{\alpha \kappa } s)}) > ct } 
\label{} \\
=& c^{-\alpha \kappa } L(p,0;ct) . 
\label{}
\end{align}
Since $ \{ X(p,0),L(p,0) \} $ is identical in law to the pair $ \{ X,L \} $ 
of the coordinate process $ X $ and its local time of $ 0 $ under $ \bP_0 $, 
we obtain {\bf (S2)}. 

Conversely, suppose that {\bf (S2)} is satisfied. 
Let us write $ \tilde{X} = \Psi_{\alpha } X $ and $ \tilde{L} = \Psi_{\alpha \kappa } L $. 
It is then obvious that for any $ t>0 $ 
\begin{align}
\int_0^t 1_{\{ \tilde{X}(s) \neq 0 \}} \d \tilde{L}(s) 
= c^{-\alpha \kappa } \int_0^{ct} 1_{\{ X(s) \neq 0 \}} \d L(s) 
= 0 . 
\label{}
\end{align}
From this it follows that $ \tilde{L} $ is a choice of the local time of $ 0 $ for $ \tilde{X} $, 
and hence $ \{ \tilde{X},\tilde{L} \} \law \{ X,k L \} $ for some constant $ k $. 
Since $ \tilde{L} \law L $, we obtain $ k=1 $. 
We denote $ \tilde{\eta} $ for the right-continuous inverse of $ \tilde{L} $ 
and define 
\begin{align}
\tilde{p}^{(l)}(t) = 
\begin{cases}
\tilde{X}(\tilde{\eta}(l-)+t) 
& \text{if $ 0 \le t < \tilde{\eta}(l)-\tilde{\eta}(l-) $}, \\
o 
& \text{if $ t \ge \tilde{\eta}(l)-\tilde{\eta}(l-) $}. 
\end{cases}
\label{}
\end{align}
We then see that the point process $ \tilde{p} = (\tilde{p}^{(l)})_{l \ge 0} $ 
is identical in law to the point process of excursions for $ X $, 
which shows that the characteristic measure of $ \tilde{p} $ is $ \vn $. 
By Lemma \ref{scaling}, 
we see that the characteristic measure of $ \tilde{p} $ is 
equal to $ c^{\alpha \kappa } \vn \circ \Psi_{\alpha }^{-1} $. 
We thus obtain {\bf (S2)$ ' $}. 

The proof is now complete. 
}

\subsection{The jumping-in vanishing case}

We prove Theorem \ref{hom}. 

\Proof[Proof of Theorem \ref{hom}]{
Set 
\begin{align}
(p^{(n)}_{\rho,j})^{(l)} = \Psi_{\alpha }^n p_{\rho,j}^{(c^{\alpha \kappa n} l)} . 
\label{}
\end{align}
By Lemma \ref{scaling}, we see that 
the characteristic measure of $ p^{(n)}_{\rho,j} $ is $ \vn^{(n)}_{\rho,j} $ 
and that 
\begin{align}
X^{(n)}_{\rho,j,\varsigma} = X(p^{(n)}_{\rho,j},\varsigma^{(n)}) 
, \quad 
L^{(n)}_{\rho,j,\varsigma} = L(p^{(n)}_{\rho,j},\varsigma^{(n)}) 
, \quad 
\eta^{(n)}_{\rho,j,\varsigma} = \eta(p^{(n)}_{\rho,j},\varsigma^{(n)}) . 
\label{}
\end{align}
Set 
\begin{align}
(\vn_n,\varsigma_n) = (\vn^{(n)}_{\rho,j},\varsigma^{(n)}) 
\quad \text{and} \quad 
(\vn_{\infty },\varsigma_{\infty }) = (\vn_{\rho^*,0},0) 
\label{}
\end{align}
and we would like to verify that 
all the assumptions of Theorem \ref{glt} are satisfied. 
It is obvious that {\bf (A1)}, {\bf (A3)} and {\bf (A4)} are satisfied. 
We have only to prove that {\bf (A2)} is satisfied. 

By the definition \eqref{eq: vnthetaj}, we have 
\begin{align}
\vn^{(n)}_{\rho,j} 
= c^{\alpha \kappa n} \vn_{\rho,0} \circ (\Psi_{\alpha }^n)^{-1} 
+ c^{\alpha \kappa n} \int_{S \setminus \{ 0 \}} j(\d x) \bP^0_x \circ (\Psi_{\alpha }^n)^{-1} . 
\label{}
\end{align}
Using \eqref{eq: vnthetaj}, {\bf (S2)$ ' $} and {\bf (S1)}, we have 
\begin{align}
\vn^{(n)}_{\rho,j} = 
\vn_{\rho,0} 
+ c^{\alpha \kappa n} \int_{S \setminus \{ 0 \}} j(\d x) \bP^0_{x_n} , 
\label{}
\end{align}
where we write $ x_n = c^{-\alpha n} x $. 
By \eqref{eq: psisigma} and {\bf (C2)}, we have 
\begin{align}
\vn_{\psi(x_n)}(T_{x_n}<T_0) 
= \sigma_{\psi(x)}(x_n) = c^{\alpha \kappa n} \sigma_{\psi(x)}(x) > 0 
\label{}
\end{align}
for $ j $-a.e. $ x \in S \setminus \{ 0 \} $. 
Using Corollary \ref{vnPo} and formula \eqref{eq: psisigma}, 
we have, for $ A \in \cB(D) $, 
\begin{align}
\vn^{(n)}_{\rho,j}(A) 
=& \vn_{\rho,0}(A) 
+ c^{\alpha \kappa n} \int_{S \setminus \{ 0 \}} \frac{j(\d x)}{\sigma_{\psi(x_n)}(x_n)} 
\cdot \vn_{\psi(x_n)} \rbra{ \{ T_{x_n} < T_0 \} \cap \theta_{T_{x_n}}^{-1} A } 
\label{} \\
=& \vn_{\rho,0}(A) 
+ \int_{S \setminus \{ 0 \}} \frac{j(\d x)}{\sigma_{\psi(x)}(x)} 
\vn_{\psi(x)} \rbra{ \{ T_{x_n} < T_0 \} \cap \theta_{T_{x_n}}^{-1} A } . 
\label{}
\end{align}
Let $ E = S \times D $. 
We define a measure $ \nu $ on $ E $ by 
\begin{align}
\nu(\d x \d w) = \delta_0(\d x) \vn_{\rho,0}(\d w) 
+ 1_{S \setminus \{ 0 \}}(x) \frac{j(\d x)}{\sigma_{\psi(x)}(x)} \vn_{\psi(x)}(\d w) . 
\label{}
\end{align}
For $ n \in \bN $, we define $ \Phi_n:E \to D $ by 
\begin{align}
\Phi_n(x,w) = 
\begin{cases}
w & \text{if $ x=0 $}, \\
\theta_{T_{x_n}}(w) & \text{if $ T_{x_n}(w) < T_0(w) $}, \\
o & \text{otherwise}. 
\end{cases}
\label{}
\end{align}
For $ n=\infty $, we define $ \Phi_{\infty }:E \to D $ by 
$ \Phi_{\infty }(x,w) = w $. 
It is obvious by {\bf (C4)} that 
{\bf (G1)} of Definition \ref{def: conv} is satisfied. 
By the definitions of $ \Phi_n $'s and $ \Phi_{\infty } $, 
we have $ T_0(\Phi_n(x,w)) \le T_0(w) $, 
and hence we see that {\bf (G5)} is satisfied. 
Since $ \| \Phi_n(x,w) \| \le \| w \| $ in any case, 
we see that {\bf (G3)$ ' $} is satisfied.

Let us verify that {\bf (G2)} and {\bf (G4)} are satisfied. 
We deal only with $ \nu $-a.e. $ (x,w) \in E $. 
\\ (i) 
For $ x = 0 $, we have, for all $ n \in \bN $, 
\begin{align}
\Phi_n(0,w) = w = \Phi_{\infty }(0,w) , 
\label{}
\end{align}
so that we have 
\begin{align}
T_0(\Phi_n(0,w)) = T_0(w) = T_0(\Phi_{\infty }(0,w)) . 
\label{}
\end{align}
(ii) 
For $ x \neq 0 $, we have 
\begin{align}
T_0(\Phi_n(x,w)) = T_0(w) - T_{x_n}(w) 
\to T_0(w) = T_0(\Phi_{\infty }(x,w)) 
\quad \text{$ \vn_{\psi(x)} $-a.e.} 
\label{}
\end{align}
For $ n \in \bN $, we define a transformation $ I_n:[0,\infty ) \to [0,\infty ) $ by 
\begin{align}
I_n(t) = 
\begin{cases}
\{ 1 + nT_{x_n}(w) \} t & \text{if $ 0 \le t < 1/n $}, \\
t + T_{x_n}(w) & \text{if $ t \ge 1/n $}. 
\end{cases}
\label{}
\end{align}
Then we easily see that $ I_n \to I $ uc. 
Since $ \Phi_n(x,w)(t) = \Phi_{\infty }(x,w)(I_n(t)) $ for $ t \ge 1/n $ 
and since $ w(0)=0 $, we obtain 
\begin{align}
\sup_{t \ge 0} \absol{ \Phi_n(x,w)(t) - \Phi_{\infty }(x,w)(I_n(t)) \Big. } 
\le 2 \sup_{0 \le t \le 1/n + T_{x_n}(w)} | w(t) | 
\to 0 . 
\label{}
\end{align}
This shows that $ \Phi_n(x,w) \to \Phi_{\infty }(x,w) $ in $ D $. 
We therefore obtain that {\bf (G2)} and {\bf (G4)} are satisfied. 

The proof is now complete. 
}

\subsection{The jumping-in dominant case}

We prove Theorem \ref{hom2}. 

\Proof[Proof of Theorem \ref{hom2}]{
Set 
\begin{align}
(p^{(n)}_{\rho,j})^{(l)} = \Psi_{\alpha }^n p_{\rho,j}^{(c^{\alpha \beta n} l)} . 
\label{}
\end{align}
By Lemma \ref{scaling}, we see that 
the characteristic measure of $ p^{(n)}_{\rho,j} $ coincides with $ \vn^{(n)}_{\rho,j} $ 
and that 
\begin{align}
X^{(n)}_{\rho,j,\varsigma} = X(p^{(n)}_{\rho,j},\varsigma^{(n)}) 
, \quad 
L^{(n)}_{\rho,j,\varsigma} = L(p^{(n)}_{\rho,j},\varsigma^{(n)}) 
, \quad 
\eta^{(n)}_{\rho,j,\varsigma} = \eta(p^{(n)}_{\rho,j},\varsigma^{(n)}) . 
\label{}
\end{align}
Set 
\begin{align}
(\vn_n,\varsigma_n) = (\vn^{(n)}_{\rho,j},\varsigma^{(n)}) 
\quad \text{and} \quad 
(\vn_{\infty },\varsigma_{\infty }) = (\vn_{0,j^*},0) 
\label{}
\end{align}
and we would like to verify that 
all the assumptions of Theorem \ref{glt} are satisfied. 
It is obvious that {\bf (A1)}, {\bf (A3)} and {\bf (A4)} are satisfied. 
We have only to prove that {\bf (A2)} is satisfied. 

Using \eqref{eq: vnthetaj}, {\bf (S2)$ ' $} and {\bf (S1)}, we have 
\begin{align}
\vn^{(n)}_{\rho,j} 
=& c^{\alpha \beta n} \vn_{\rho,0} \circ (\Psi_{\alpha }^n)^{-1} 
+ c^{\alpha \beta n} \int_{S \setminus \{ 0 \}} j(\d x) \bP^0_x \circ (\Psi_{\alpha }^n)^{-1} 
\label{} \\
=& c^{-\alpha (\kappa - \beta) n} \vn_{\rho,0} 
+ c^{\alpha \beta n} \int_{S \setminus \{ 0 \}} j(\d x) \bP^0_{x_n} , 
\label{}
\end{align}
where we write $ x_n = c^{-\alpha n} x $. 
Using \eqref{eq: psisigma}, {\bf (C2)}, 
Corollary \ref{vnPo}, formula \eqref{eq: psisigma} and {\bf (C5)}, 
we have, for $ A \in \cB(D) $, 
\begin{align}
\vn^{(n)}_{\rho,j}(A) 
=& c^{-\alpha (\kappa - \beta) n} \vn_{\rho,0}(A) 
+ c^{\alpha \beta n} \int_{S \setminus \{ 0 \}} 
\frac{j(\d x)}{\sigma_{\psi(x)}(x_n)} 
\vn_{\psi(x)} \rbra{ \{ T_{x_n} < T_0 \} \cap \theta_{T_{x_n}}^{-1} A } 
\label{} \\
=& c^{-\alpha (\kappa - \beta) n} \vn_{\rho,0}(A) 
+ \int_{S''} \mu(\d y) 
\vn_{\psi(J^*y)} \rbra{ \{ T_{J_ny} < T_0 \} \cap \theta_{T_{J_ny}}^{-1} A } . 
\label{}
\end{align}
Let $ E = ((0,\infty ) \cup S'') \times D $. 
We define a measure $ \nu $ on $ E $ by 
\begin{align}
\nu(\d y \d w) = 1_{(0,\infty )}(y) \d y \vn_{\rho,0}(\d w) 
+ 1_{S''}(y) \mu(\d y) \vn_{\psi(J^*y)}(\d w) . 
\label{}
\end{align}
For $ n \in \bN $, we define $ \Phi_n:E \to D $ by 
\begin{align}
\Phi_n(y,w) = 
\begin{cases}
w & \text{if $ y \in (0,c^{-\alpha (\kappa - \beta) n}) $}, \\
\theta_{T_{J_ny}}(w) 
& \text{if $ y \in S'' $ and $ T_{J_ny}(w)<T_0(w) $}, \\
o & \text{otherwise}. 
\end{cases}
\label{}
\end{align}
For $ n=\infty $, we define $ \Phi_{\infty }:E \to D $ by 
\begin{align}
\Phi_{\infty }(y,w) = 
\begin{cases}
\theta_{T_{J^* y}}(w) & \text{if $ y \in S'' $ and $ T_{J^*y}(w)<T_0(w) $}, \\
o & \text{otherwise}. 
\end{cases}
\label{}
\end{align}
By the above argument and by {\bf (C6)}, 
we find that 
{\bf (G1)} of Definition \ref{def: conv} is satisfied. 
By the definitions of $ \Phi_n $'s and $ \Phi_{\infty } $, 
we have $ T_0(\Phi_n(y,w)) \le T_0(w) $, 
and hence we see that {\bf (G5)} is satisfied. 
Since $ \| \Phi_n(y,w) \| \le \| w \| $ in any case, 
we see that {\bf (G3)$ ' $} is satisfied. 

Let us verify that {\bf (G2)} and {\bf (G4)} are satisfied. 
We deal only with $ \nu $-a.e. $ (y,w) \in E $. 
\\ (i) 
For $ y \in (0,\infty ) $, we have 
\begin{align}
\Phi_n(y,w) = 
\begin{cases}
w & \text{if $ n < \frac{1}{\alpha (\kappa -\beta ) \log c} \log \frac{1}{y} $}, \\
o & \text{otherwise}, 
\end{cases}
\label{}
\end{align}
since $ c>1 $. 
This shows that $ \Phi_n(y,w)=o=\Phi_{\infty }(y,w) $ for large $ n $. 
\\ (ii) 
For $ y \in S'' $, we have $ T_{J_ny}(w) \to T_{J^*y}(w) $ by {\bf (C5)}. 
We define transformations $ I_n:[0,\infty ) \to [0,\infty ) $ for $ n \in \bN $ as follows. 
Let 
\begin{align}
\tau_n(w) = T_{J_ny}(w) - T_{J^*y}(w) 
\quad \text{and} \quad 
\tau^{\pm}_n(w) = \max \{ \pm \tau_n(w),0 \} . 
\label{}
\end{align}
For $ n \in \bN $, we define 
\begin{align}
I_n(t) = 
\begin{cases}
\frac{1/n + \tau^+_n(w)}{1/n + \tau^-_n(w)} t & \text{if $ 0 \le t < 1/n + \tau^-_n(w) $}, \\
t + \tau_n(w) & \text{if $ t \ge 1/n + \tau^-_n(w) $}. 
\end{cases}
\label{}
\end{align}
Then we easily see that $ I_n \to I $ uc. 
Let us write 
\begin{align}
T^n_{\rm min}(w) = \min \{ T_{J_ny}(w),T_{J^*y}(w) \} 
, \quad 
T^n_{\rm max}(w) = \max \{ T_{J_ny}(w),T_{J^*y}(w) \} . 
\label{}
\end{align}
Since $ \Phi_n(y,w)(t) = \Phi_{\infty }(y,w)(I_n(t)) $ for $ t \ge 1/n + \tau^-_n(w) $, 
we obtain 
\begin{align}
\sup_{t \ge 0} \absol{ \Phi_n(y,w)(t) - \Phi_{\infty }(y,w)(I_n(t)) \Big. } 
\le 2 \sup_{T^n_{\rm min}(w) \le t \le 1/n + T^n_{\rm min}(w) } | w(t) | . 
\label{}
\end{align}
The last quantity converges to 0 by {\bf (C5)}. 
This shows that $ \Phi_n(y,w) \to \Phi_{\infty }(y,w) $ in $ D $. 
We therefore obtain that {\bf (G2)} and {\bf (G4)} are satisfied. 

The proof is now complete. 
}

\def\cprime{$'$}

\end{document}